\documentclass[%
 aip,
 amsmath,amssymb,
 reprint,%
]{revtex4-1}

\usepackage{graphicx}
\usepackage{dcolumn}
\usepackage{bm}

\usepackage[utf8]{inputenc}
\usepackage[T1]{fontenc}
\usepackage{mathptmx}
\usepackage{etoolbox}
\usepackage[colorlinks, linkcolor=red, anchorcolor=blue, citecolor=green]{hyperref}
\usepackage{url}

\usepackage{amsmath}
\usepackage{booktabs}  
\usepackage{graphicx}
\usepackage{subfigure}
\usepackage{framed,multirow}
\usepackage{latexsym}
\usepackage[ruled,vlined]{algorithm2e}
\usepackage{bm}
\makeatletter
\def\@email#1#2{%
 \endgroup
 \patchcmd{\titleblock@produce}
  {\frontmatter@RRAPformat}
  {\frontmatter@RRAPformat{\produce@RRAP{*#1\href{mailto:#2}{#2}}}\frontmatter@RRAPformat}
  {}{}
}%
\makeatother
\begin{document}

\preprint{AIP/123-QED}

\title[Sample title]{Critical grid method: An extensible Smoothed Particle Hydrodynamics fluid general interpolation method for Fluid-Structure Interaction surface coupling based on preCICE}
\author{Sifan Long}
 \affiliation{%
College of Computer, Central South University, 410083, ChangSha, China.}%
 \altaffiliation[Also at ]{College of Computer, National University of Defense Technology, 410073, ChangSha, China.}
\author{Xiaowei Guo}%
\affiliation{Institute for Quantum Information $\&$ State Key Laboratory of High Performance Computing, College of Computer, National University of Defense Technology, 410073, ChangSha, China.}%

\author{Xiaokang Fan}
 \email{fanxiaokang@nudt.edu.cn}
 \homepage{http://www.Second.institution.edu/~Charlie.Author.}
\affiliation{%
College of Computer, National University of Defense Technology, 410073, ChangSha, China.}%

\author{Canqun Yang}
\affiliation{%
National Supercomputing Center in Tianjin, 300457, TianJin, China.}%

\date{\today}

\begin{abstract}
Solving Fluid-Structure Interaction (FSI) problems using traditional methods is a big challenge in the field of numerical simulation. As a powerful multi-physical field coupled library, preCICE has a bright application prospect for solving FSI, which supports many open/closed source software and commercial CFD solvers to solve FSI problems in the form of a black box. However, this library currently only supports mesh-based coupling schemes. This paper proposes a critical grid (mesh) as an intermediate medium for the particle method to connect a bidirectional coupling tool named preCICE. The particle and critical mesh are used to interpolate the displacement and force so that the pure Lagrangian Smoothed Particle Hydrodynamic (SPH) method can also solve the FSI problem. This method is called the particle mesh coupling (PMC) method, which theoretically solves the mesh mismatch problem based on the particle method to connect preCICE. In addition, we conduct experiments to verify the performance of the PMC method, in which the fluid and the structure is discretized by SPH and the Finite Element Method (FEM), respectively. The results show that the PMC method given in this paper is effective for solving FSI problems. Finally, our source code for the SPH fluid adapter is open-source and available on GitHub\footnote{\href{https://github.com/terrylongsifan/AdapterSPH}{www.github.com/terrylongsifan/AdapterSPH}} for further developing preCICE compatibility with more meshless methods.
\end{abstract}

\maketitle

\section{Introduction}\label{section1}
In the field of engineering applications, multi-physical field problems involve many research objects, including acoustics, electricity, mechanics, optics, thermology and their cross fields. Fluid-Structure Interaction (FSI) is a common phenomenon in nature \cite{griffith2020immersed}, which is an important content in the study of multi-physical fields. Due to the complex boundary conditions and changes in the flow field, its solution presents highly nonlinear characteristics \cite{liu2019smoothed}, so that it is impossible to directly solve its theoretical value. Therefore, the simulation of various complex FSI scenarios is facing great challenges in the academia and industry, Fluid-Structure Interaction has become a very hot research topic in the field of multi-physical coupling \cite{kim2019immersed,mehryan2020fluid,toma2021fluid}. At present, with the development of computer hardware equipment and the improvement of computing performance, Computational Fluid Dynamics (CFD) and Computational Solid Mechanics (CSM) have achieved rapid development in recent decades \cite{shen2020recent,jain2020advances}, especially various high-performance CFD and CSM algorithms have been proposing and continuously improving, which provides strong support for computer-aided engineering (CAE) to solve fluid-structure interaction problems.

At present, mesh-based numerical methods are mostly used to solve FSI problems, among which the representative methods are Finite Element Method (FEM) \cite{jagota2013finite,szabo2021finite}, Finite Volume Method (FVM) \cite{moukalled2016finite,lin2013finite} and Finite Difference Method (FDM) \cite{patidar2016nonstandard}. They all use similar methods to discretize the governing equations of fluid flow. First, the solution domain is divided into non-overlapping meshs or elements, and then the governing equations of fluid are discretized on the meshs or elements. Therefore, the quality of the meshs or elements determines the accuracy of the calculation results. This method has achieved great success in solving FSI problems. For example, Fong et al. used finite element analysis in oil and gas industry to solve the FSI problem of damping flow pipeline vibration \cite{fong2017fluid}. In addition, Pedro and Qian\cite{martinez2018efficient} used the Finite Volume Method to solve the strong coupling problem between two-phase flow and isothermal non-linear elastic bodies under the FSI framework, and they successfully analyzed the interaction mechanism between fluid and elastic structure under the benchmark example, which fully verified the potential of this method in engineering application.
At the same time, the use of a variety of mesh-based methods has also become a popular solution for solving FSI problems. Peterson \textit{et al.}\cite{peterson2020strongly} has solved the lubrication problem of two-dimensional elastohydrodynamics by combining finite element and finite volume methods. This solution idea can be compatible with the advantages of many methods. Therefore, the mesh-based FSI solution methods have been widely used in mechanical engineering, aerospace, ocean engineering, biomedicine, energy development and other fields \cite{hou2012numerical}.

Although the mesh-based method has achieved great success in the application of FSI. However, with the in-depth development of the theoretical and the broadening of engineering applications, the mesh-based method also shows many drawbacks that are difficult to overcome \cite{nguyen2008meshless}.
They have encountered great challenges in practical applications, which has seriously hindered their development. The generation and management of mesh have become a complex and arduous task in the solution process. At the same time, the requirements of mesh generation without overlap, distortion, and entanglement must be met \cite{thompson2020structured}. The quality of the mesh directly determines the final simulation results.
In particular, when dealing with the hot issues of moving boundary and interface discontinuity concerned by current engineering applications, such as large deformation of fluid, high-speed impact response of materials, explosion and impact, the mesh-based method is easy to cause mesh distorted in the calculation procedure \cite{shang2017numerical}, resulting in negative density of grid elements, which leads to very small time step and termination of calculation.

In recent years, in order to overcome the inherent defects of mesh based methods, various meshless methods based on node interpolation have been proposed. These methods with Lagrangian characteristics do not rely on the background grid, so they can easily deal with the problems that are difficult to deal with based on mesh-based methods, such as discontinuous interfaces, complex geometric boundaries, incompressible flows with free surfaces. In addition, it also has the advantages of high continuity of approximate functions, easy to add or delete nodes and easy to implement adaptive technology, which makes the mesh free method quickly become a hot spot in the research of numerical methods. At the same time, many mesh free methods have emerged, such as Smooth Particle Hydrodynamics (SPH) method \cite{gingold1977smoothed}, Diffuse Element Method (DEM) \cite{1992Generalizing} and Finite Point Method (FPM) \cite{onate1996finite}. Among them, as a typical representative of mesh free method, the Smooth Particle Hydrodynamics method was proposed by Monaghan \textit{et al.}\cite{gingold1977smoothed} and successfully applied to solve three-dimensional open boundary astrophysical problems in 1977. 
Its principle is to use a series of node particles with physical information such as position, mass, and velocity in the solution domain to discrete the field variables in the space, and then solve the evolution process of the whole simulation results. Unlike other meshless methods, which only use node particles as interpolation nodes in the procedure of calculation, the node particles of the SPH method also carry the information of materials \cite{liu2010smoothed}.
Therefore, due to these characteristics, the SPH method is very suitable for simulating large deformation of materials \cite{rahman2021different}, impact fracture of solid media \cite{meng2021advances}, multiphase flow \cite{ye2019smoothed}, and Fluid-Structure Interaction processes. This method has attracted the extensive attention of researchers all over the world. It has been widely used in the fields of energy mining, biomedicine and aerospace \cite{liu2010smoothed}. The relevant SPH community codes can be download from (\href{https://www.spheric-sph.org/}{www.spheric-sph.org}).
Although SPH method has been successfully applied in the field of numerical simulation. However, it still has many defects, such as low computational efficiency, tensile instability, and difficult to deal with boundary conditions. These shortcomings make SPH still face great challenges in the application of Fluid-Structure Interaction. Therefore, coupling SPH method with other mature mesh-based methods is an important research content that SPH can be more widely used to solve FSI problems.

Previously, we introduced some theoretical characteristics of mesh-based methods and meshless methods, including their advantages, disadvantages and their respective application ranges. For solving the problem of FSI, the force of fluid on the structure causes the deformation of the structure part, and the deformation of the structure will affect the flow of fluid. This is a highly nonlinear system with two-way coupling, which is difficult to be solved directly by using the mesh-based method or the meshless method alone. Due to the complementary feature of the Lagrangian grid method and Eulerian grid method, based on this feature, the coupled Eulerian-Lagrangian (CEL) \cite{ireland2017improving} method and Arbitrary Lagrangian-Eulerian (ALE) \cite{barlow2016arbitrary} method are proposed. They take into account both the advantages of the two methods and can avoid their respective disadvantages. However, although CEL and ALE methods are used, the simulation results still produce errors when the mesh is highly distorted. At present, many studies\cite{long2017arbitrary,fragassa2019dealing,yang2012free} have shown that mesh method combined with particle method usually has many advantages in solving FSI problems, which has become a promising direction for solving FSI problems. On the one hand, the Finite Element Method is usually used to solve the structural module of FSI because of its high accuracy, high computational efficiency, and strong robustness. On the other hand, SPH method is usually used to solve the fluid module of FSI because of its excellent ability to deal with large deformation of materials and free surface flow. According to practical problems, mesh elements and particle nodes are solved iteratively through consolidation, contact, and transformation algorithms.
Therefore, after being compatible with the advantages of the mesh method and the particle method, the SPH-FEM coupling framework has achieved rapid development.  Wang and  Wu\cite{wang2021numerical} established a dynamic compaction model through the FEM-SPH coupling method. Their tests showed that the effect is better than using the traditional FEM or pure SPH.
Fragassa \textit{et al.}\cite{fragassa2019dealing} used the FEM-SPH coupling framework to realize the multi field coupling between water flow, air and structure. 
In addition, Jayasinghe \textit{et al.}\cite{2020Impact} used the FEM-SPH coupling model to study the punching of piles under lateral load and the impact on adjacent piles, and the model can avoid dealing with the numerical oscillation caused by mesh tangling and remeshing, which improved the simulation accuracy and calculation efficiency.

Although the mesh-based method coupled particle method can be successfully applied to FSI problems, most of the current coupling methods rely heavily on the problem background to be solved, and the code is difficult to be directly applied to solve the new FSI problems. In particular, in the face of different materials and different coupling scenarios (such as contact and penetration), it requires strong skills to deal with the contact algorithm between mesh elements and particle nodes. Therefore, the framework of program is usually complex and the reusability is poor. In order to overcome these shortcomings, the Precise Code Interaction Coupling Environment (preCICE) provides a black box coupling library to solve multiple physical field problems \cite{bungartz2016precice}, it allows non-coupling-experts to spend the least effort to achieve a high-precision and stable coupling solution, avoiding the repeated implementation of the coupling algorithm. This software is an open-source surface coupling library written by C++ to connect different solvers. Because of its rich documentation and powerful functions. It can support various types of meshes (Eulerian, Lagrangian, structural mesh, unstructured mesh, static and dynamic mesh), and also support a variety of discrete methods (Finite Element, Finite Volume and Finite Difference). Different from the previous client-server coupling mode, preCICE directly supports coupling communication between multiple solvers. It not only successfully solves many problems of FSI, but also it has been widely used in the fields of heat transfer \cite{rodenberg2021fenics}, acoustics \cite{blom2016partitioned}, and other multiple physical fields. Therefore, it is widely spread in the field of multiple physical fields. The fluid/solid solver and preCICE are connected by specific adapters. This highly modular design enables preCICE to easily use the coupling scheme without major changes to the core code of the solver. Therefore, preCICE couples a large number of influential solver software in the field of numerical simulation. Including fluid solver OpenFOAM \cite{jasak2009openfoam} and SU2 \cite{economon2016su2}, it also includes the solid solver Deal.II \cite{bangerth2007deal}, FEniCS \cite{alnaes2015fenics}, and commercial software ANSYS Fluent \cite{pashchenko2018ansys} and COMSOL \cite{dickinson2014comsol}. At the same time, more and more solvers provide methods compatible with preCICE coupling framework. Relevant software and corresponding adapters can be obtained on (\href{http://precice.org/}{www.precice.org}).

preCICE provides a convenient coupling scheme with its perfect functions and user-friendly development documents, making it one of the open-source software that has attracted much attention in the field of multi physical fields. However, the current research on docking multiple physics coupling framework preCICE based on meshless method (eg. SPH) needs to be further expanded, which makes the application potential of meshless method that is good at dealing with large deformation in FSI problems be effectively explored. It is an important research content to develop an adaptive algorithm that is compatible with the meshless method and supports high-precision Fluid-Structure Interaction algorithm in the preCICE coupling framework.
Therefore, this paper proposes a method for docking the SPH particle adaptation mesh of preCICE coupling framework, which is called the particle-mesh coupling (PMC) algorithm. Based on this interpolation method, we successfully realize the coupling of the SPH method and Finite Element Method, and it is applied to the solution of FSI problems. Because the coupling mechanism is universal, it can be coupled with any solver compatible with preCICE in theory. It provides a solution that uses particle (or meshless) methods in the influential open-source preCICE. This scheme can quickly make developer without coupling background knowledge become coupling experts, and let them focus on the problems of their own industries instead of spend lot of time for writting code.

The remaining structure of this paper is mainly divided into the following sections: Firstly, Section \ref{section2} introduces the basic principle of the Smoothed Particle Hydrodynamics method, and Section \ref{section3} introduces the working mechanism and current development of preCICE coupling framework. Then, in Section \ref{section4}, a method based on Particle-Mesh Coupling (PMC) is proposed, which couples the Particle-mesh method to the preCICE framework. Section \ref{section5} describes the specific implementation procedure of the PMC method in detail. In Section \ref{section6}, a numerical example is given to test the proposed algorithm and verify the performance of the PMC method. Finally, Section \ref{section7} is the conclusion and the future outlook.

\section{Smoothed particle hydrodynamics method}\label{section2}
In this section, firstly, we will introduce the interpolation theory of the SPH method, including kernel approximation and particle approximation. Then, the procedure of discretization of the conservation equation controlling fluid flow using the SPH method is described in detail. Finally, the discretized approximate solution is solved by numerical integration algorithm.
\subsection{Kernel approximation}
In the SPH method, for any continuous field function $f(\bm{r})$, the function value at a point on the definition domain $\Omega$ can be expressed by the following equation (\ref{basic_eq1}).
\begin{equation}\label{basic_eq1}
f(\bm{r})=\int_{\Omega }f(\bm{r}')\delta (\bm{r}-\bm{r}')d\bm{r}'
\end{equation}

Where $\bm{r}$ is the vector in $n$-dimensional space ($n=1,2,3,\cdots $) and $\delta (\bm{r}-\bm{r}')$ is Dirac $\delta$ Function, which satisfies the following properties.
\begin{equation}\label{basic_eq2}
\delta \left ( \bm{r}-\bm{r}' \right ) =\left\{\begin{matrix}
\infty, \bm{r}=\bm{r}' 
\\
0, \bm{r}\ne \bm{r}'
\end{matrix}\right.
\end{equation}

\begin{equation}\label{basic_eq3}
\int_{\mathbb{R}^{n}}\delta\left ( \bm{x}\right )d\bm{x}=1
\end{equation}

However, Dirac $\delta$ function is discontinuous in reality, so the smooth function $W\left (\bm{r}-\bm{r}',h \right )$ is used to approximately replace the $\delta$ function, which obtains the approximate kernel function expression of the field function $f(\bm{r})$.
\begin{equation}\label{basic_eq4}
f(\bm{r})\approx \int_{\Omega }f({\bm{r}}')W(\bm{r}-{\bm{r}}',h)d{\bm{r}}'
\end{equation}

$W(\bm{r}-{\bm{r}}',h)$ represents the smooth kernel function in SPH method, and its value depends on the distance $\left | \bm{r}-{\bm{r}}' \right |$ for two points and the smooth length $h$, and its influence region is determined by the smoothing factor $\kappa$ and the smoothing length $h$.
\begin{figure}[!t]
	\centering
	\includegraphics[scale=.8]{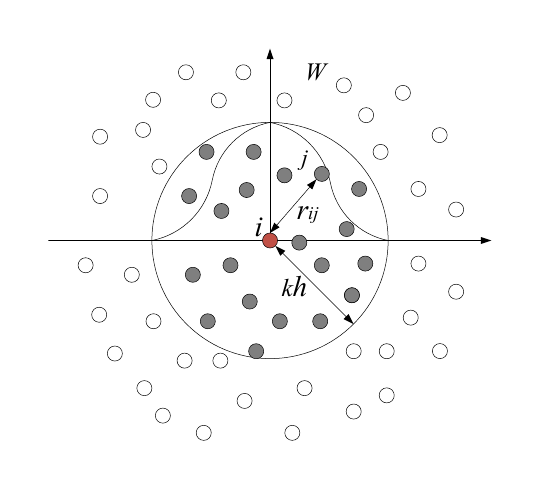}
	\caption{Kernel function approximation and particle approximation principle of the SPH method.}
	\label{fig1}
\end{figure}

Based on the same principle, the kernel approximation of the derivative function of $f(\bm{r})$ can be evaluated according to the integration formula and the divergence theorem.
\begin{equation}\label{basic_eq5}
\nabla f(\bm{r})\approx-\int_{\Omega }f(\bm{r}') \cdot \nabla W(\bm{r}-\bm{r}',h)d\bm{r}'
\end{equation}

Equations (\ref{basic_eq4}) and (\ref{basic_eq5}) show that any field function and its derivatives in space can be expressed by the smooth function $W$. However, equation (\ref{basic_eq5}) can only show that the support region completely falls within the boundary range. When the support region is truncated by the boundary, boundary conditions need to be imposed. Therefore, imposing boundary conditions is also a hot and difficult topic in the theoretical research of the SPH method.

\subsection{Particle approximation}
Equations (\ref{basic_eq4}) and (\ref{basic_eq5}) obtain the integral expressions of the field function and its derivatives based on the smooth kernel function $W$. But it can not discretize the governing equations of the fluid. Therefore, the SPH method uses the particle approximation method to interpolate and discretize the smooth kernel function. Its discretization principle is as follows.

In the solution domain $\Omega$, as shown in Fig.\ref{fig1}, the integral expression of the kernel function can be further transformed into the form of weighted summation of a series of discrete particles in the support domain, which is called the particle approximation of the SPH method. If $V_j$ is used to represent the volume element $d{\bm{r}}'$ of particle $j$ with arbitrary distribution, the mass of particle $j$ can be expressed by equation (\ref{basic_eq6}).
\begin{equation}\label{basic_eq6}
m_{j}=\Delta V_{j}\rho _{j}
\end{equation}

Where $\rho_j$ is the density of particle $j$. Based on equations (\ref{basic_eq4}) and (\ref{basic_eq6}), the particle estimation expression for the field function at particle $i$ can be solved as shown in equation (\ref{basic_eq7}).
\begin{equation}\label{basic_eq7}
f(\bm{r}_{i})\approx \sum_{j}f(\bm{r}_{j})\frac{m_{j}}{\rho_{j}}W(\bm{r}_{i}-\bm{r}_{j},h)
\end{equation}

Similarly, according to equations (\ref{basic_eq5}) and (\ref{basic_eq6}), the derivative expression of the field function at particle $i$ can be obtained as shown in equation (\ref{basic_eq8} or \ref{basic_eq11}).
\begin{equation}\label{basic_eq8}
\nabla \cdot f(\bm{r}_i)\approx \rho_i {\sum_{j=1}^{N}} m_j\left [ \frac{f(\bm{r}_j)}{\rho_{j}^{2}} + \frac{f(\bm{r}_i)}{\rho_{i}^{2}} \right ]\cdot \nabla_iW_{ij}
\end{equation}

\begin{equation}\label{basic_eq11}
\nabla \cdot f(\bm{r}_i)\approx \frac{1}{\rho_i} \sum_{j=1}^{N} m_j\left [ f(\bm{r}_j)-f(\bm{r}_i) \right]\cdot \nabla_iW_{ij}
\end{equation}

Where 
\begin{equation}\label{basic_eq9}
\nabla_i W_{ij} = \frac{\bm{r}_{ij}}{r_{ij}} \cdot \frac{\partial W_{ij}}{\partial r_{ij}}
\end{equation}

\begin{equation}\label{basic_eq10}
r_{ij}=\left | \bm{r}_{ij} \right | =\left | \bm{r}_i-\bm{r}_j \right | 
\end{equation}

$r_{ij}$ represents the distance between particles $i$ and $j$, then, an excellent feature of equation (\ref{basic_eq10}) is that the right side of the equation appears in the form of paired particles. The advantage of this symmetric form is that it can improve the calculation accuracy of the SPH method \cite{2006Restoring}.

\subsection{Momentum equation}
In the field of Computational Fluid Dynamics, the momentum conservation equation controlling fluid flow in a continuous field can be described by equation (\ref{basic_eq12}).
\begin{equation}\label{basic_eq12}
\frac{\mathrm{d} \bm{v}}{\mathrm{d} t} =-\frac{1}{\rho} \nabla P+\bm{g}+\bm{\Gamma}
\end{equation}

Where $\bm{\Gamma}$ is the dissipative term and $\bm{g}$ is the gravitational acceleration ($\bm{g}=9.8m/s^2$). Considering the effect of fluid viscosity and surface tension, artificial viscosity was proposed by Monaghan to simulate fluid flow because it is easy to realize. Equation (\ref{basic_eq12}) can be written as equation (\ref{basic_eq13}).
\begin{equation}\label{basic_eq13}
\frac{\mathrm{d} \bm{v}_i}{\mathrm{d} t} =-\sum_{j=1}^{N} m_j\left (\frac{P_j}{\rho_{j}^{2} } +\frac{P_i}{\rho_{i}^{2} } +\Pi_{ij} \right ) \nabla_iW_{ij}+\bm{g}
\end{equation}

Where $P$ and $\rho$ are the pressure and density of the corresponding particle $i$ or $j$, respectively, and $\Pi_{ij}$ is the viscosity term, which is given by equation (\ref{ns_eq1}).
\begin{equation}\label{ns_eq1}
\Pi_{ij}=\left\{\begin{matrix}
\frac{-\alpha \bar{c_{ij}}\mu_{ij}}{}  ,\bm{v}_{ij}\cdot \bm{r}_{ij}< 0
\\
0,\bm{v}_{ij}\cdot \bm{r}_{ij}< 0
\end{matrix}\right.
\end{equation}

Where $\bm{r}_{ij}=\bm{r}_i-\bm{r}_j$, $\bm{v}_{ij}=\bm{v}_i-\bm{v}_j$, $\mu_{ij}=h\bm{v}_{ij}\cdot \bm{r}_{ij}/(r_{ij}^2+\eta^2)$, $\bar{ c_{ij}}=0.5(c_i+c_j)$, $c$ is the speed of local sound, $h$ is the smooth length, $\eta^2=0.01h^2$, $\alpha$ is a coefficient associated with the problem.

In addition, in most application scenarios, the laminar viscous stresses also needs to be introduced into the momentum equation \cite{Edmond2002Simulation}, which can be written as 
\begin{equation}\label{ns_eq2}
\left ( \upsilon_0\nabla^2\bm{v} \right )_i=\sum_{j=1}^{N} m_j\left ( \frac{4\upsilon_0 \bm{r}_{ij}\cdot \nabla_i W_{ij}}{(\rho_i+\rho_j)(r_{ij}^2+\eta^2)}  \right ) \bm{v}_{ij}
\end{equation}

Where $\upsilon_0$ is the kinetic viscosity, generally, its value is $10^{-6}m^2s$ for water. Therefore, according to the interpolation principle of the SPH method, equation (\ref{ns_eq2}) can be transformed into equation (\ref{ns_eq3}) with discrete form of the SPH method through equations (\ref{basic_eq7}) and (\ref{basic_eq8}).
\begin{multline}\label{ns_eq3}
\frac{\mathrm{d} \bm{v}}{\mathrm{d} t}=-\sum_{j=1}^{N}m_j\left ( \frac{P_j}{\rho_j^2}+ \frac{P_i}{\rho_i^2} \right )\nabla_iW_{ij}+\bm{g}+ \\
\sum_{j=1}^{N} m_j\left ( \frac{4\upsilon_0 \bm{r}_{ij}\cdot \nabla_i W_{ij}}{(\rho_i+\rho_j)(r_{ij}^2+\eta^2)}  \right ) \bm{v}_{ij}
\end{multline}

Gotoh\cite{2001Sub} first introduced the concept of Sub-Particle Scale (SPS) to describe the affected by turbulence for their Moving Particle Semi-implicit (MPS) model. Their momentum conservation equation is defined as follows.
\begin{equation}\label{ns_eq4}
\frac{\mathrm{d} \bm{v}}{\mathrm{d} t}=-\frac{1}{\rho} \nabla P+\bm{g}+\upsilon_0 \nabla^2\bm{v}+\frac{1}{\rho}\nabla \cdot \vec{\tau}  
\end{equation}

Where equation (\ref{ns_eq2}) is used as the laminar term, and $\vec{\tau}$ is the SPS stress tensor. Finally, Rogers \textit{et al.}\cite{dalrymple2006numerical} introduce SPS into the weakly compressible SPH (WCSPH) method by using Favre averaging, and the SPH discrete form of equation (\ref{ns_eq4}) can be written as
\begin{multline}\label{ns_eq5}
\frac{\mathrm{d} \bm{v}}{\mathrm{d} t}=-\sum_{j=1}^{N}m_j\left ( \frac{P_j}{\rho_j^2}+ \frac{P_i}{\rho_i^2} \right )\nabla_iW_{ij}+\bm{g}+ \\
\sum_{j=1}^{N} m_j\left ( \frac{4\upsilon_0 \bm{r}_{ij}\cdot \nabla_i W_{ij}}{(\rho_i+\rho_j)(r_{ij}^2+\eta^2)}  \right ) \bm{v}_{ij}+ \linebreak \\
\sum_{j=1}^{N} m_j\left ( \frac{\vec{\tau}_{xy}^j }{\rho_j^2}+\frac{\vec{\tau}_{xy}^i }{\rho_i^2}  \right ) \nabla_i W_{ij}
\end{multline}

\subsection{Continuity equation}
In the process of solving weakly compressible SPH simulations, the mass of particles remains unchanged, the initial density solution is obtained by weighted summation of particle masses \cite{monaghan1992smoothed}, which tends to cause violent fluctuations near the boundary and free surface. Therefore, the density can be solved currently by the conservation equation of mass, or continuity equation as shown in equation (\ref{ns_eq6}).
\begin{equation}\label{ns_eq6}
\frac{\mathrm{d} \rho_i}{\mathrm{d} t} =\sum_{j=1}^{N} m_j \bm{v}_{ij}\cdot \nabla_i W_{ij}
\end{equation}

\subsection{Equation of state}
In the weakly compressible SPH method, the pressure is solved by the equation of state (EOS) \cite{monaghan1994simulating}. Therefore, the relationship between pressure and density is described in equation (\ref{ns_eq7}).
\begin{equation}\label{ns_eq7}
P=c_{0}^2 \rho_0/\gamma\left [ \left ( \frac{\rho}{\rho_0}  \right )^{\gamma }-1 \right ] 
\end{equation}

Where $\gamma=7$ for water, $\rho_0=1000 kg/m^3$ is the reference density, and $c_0=c(\rho_0)=\sqrt{\partial P/\partial \rho} $ is the local speed of sound when $\rho=\rho_0$.

\section{Precise code interaction coupling environment}\label{section3}
\begin{figure}[!t]
	\centering
	\includegraphics[scale=.45]{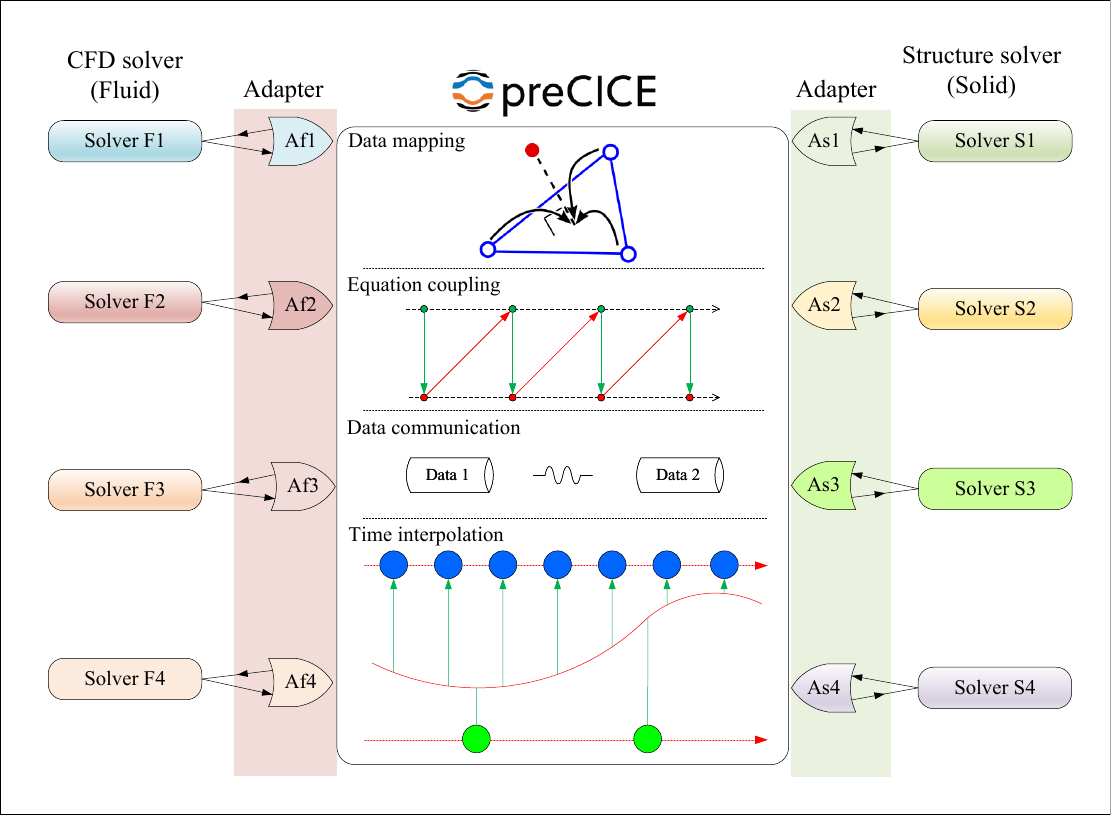}
	\caption{preCICE is an efficient coupled library for multiphysical field simulation, which provide data mapping between non-matching grids, communication, and equation coupling schemes. Solver and preCICE are connected through their adapters (eg. Afx and Asx, $x=1,2,3,\cdots$), which call the interfaces provided by preCICE to exchange data between solids and fluids.}
	\label{fig2}
\end{figure}

This section mainly introduces the theoretical principle and technical implementation of preCICE framework. As a general coupling framework, preCICE provides a fully functional coupling library that is compatible with multiple solvers to realize the simulation of multiple physical fields as shown in Fig.\ref{fig2}. 
It allows coupling solution between different solvers. preCICE provides many interfaces to connect external data for exchange, which is called by the adapter that connected with the solver. For FSI problems, existing components are often used in fluid and solid calculations. For example, the fluid is solved by SPH, and the solid is solved by FEM. Therefore, it needs to provide SPH-adapter and FEM-adapter to connect to preCICE. Then, preCICE provides a coupling scheme to complete the solution of fluid and solid parts. It includes the interpolation of non-matching grid data, coupling scheme, time step advance, data communication and other functions.

Thanks to its high degree of encapsulation and modularity, developers do not need to understand its internal details and use it as a black-box model, which greatly reduces the research threshold for realizing the coupling of multiple physical fields.

\subsection{Coupling schemes}
preCICE provides two coupling schemes to realize the coupling process: 1) explicit coupling schemes and 2) implicit coupling schemes. The difference lies in whether all physical quantities are obtained at the same time. The explicit coupling scheme needs to be called in a single time step, while the implicit coupling scheme needs to iterate the coupling equation until it converges, so it has iterative subcycles.
Taking the solution of two solvers on FSI problem as an example, the fluid solver is annotated as $F_1$, and the solid solver is denoted as $S_2$. The data of vectorization type exchanged between the two solvers are $X_1$ and $X_2$, respectively. Therefore, they have the following corresponding relationship on the coupling surface.
\begin{equation}\label{ns_eq8}
F_1:X_1\Longrightarrow X_2; S_2:X_2\Longrightarrow X_1
\end{equation}

Most of the partitioned coupling schemes belong to Dirichlet–Neumann type coupling \cite{denk2015structurally}. In FSI problem, fluid solver $F_1$ needs to provide forces $X_1$ as the input of the wet interface and transmit it to solid solver $S_2$ through preCICE. Then, solid solver $S_2$ sends back displacements $X_2$ to fluid solver $F_1$. Therefore, boundary conditions and conservation equations are also satisfied.

\subsubsection{Explicit coupling schemes}
preCICE provides an explicit coupling scheme and its parallel method. In general, the explicit coupling scheme is based on the traditional interleaving scheme, which requires solver $F_1$ to use the old boundary value $x_1^{(n)}$ in time steps from $t_n$ to $t_{n+1}$ to calculate the value $x_2^{(n+1)}$ of the next time step. 
\begin{equation}\label{ns_eq9}
x_2^{(n+1)}=F_1^{(n)}(x_1^n)
\end{equation}

Meanwhile, the solver $S_2$ waits for the calculation result of $F_1$ and updates the value of $x_2$ according to its return value and the value of $x_1$ for the next time step ($x_1^{(n+1)}$).
\begin{equation}\label{ns_eq10}
x_1^{(n+1)}=S_2^{(n)}(x_2^{(n+1)})
\end{equation}          

In an explicit parallel scheme, two solvers can simultaneously compute the value of time $t_n$. Thus, communication may occur after they have completed their calculations for reducing waiting times. The explicit parallel scheme can be described by equation (\ref{ns_eq11}).
\begin{equation}\label{ns_eq11}
\left\{\begin{matrix}
x_2^{(n+1)}=F_1^{(n)}(x_1^n)  \\
\\
x_1^{(n+1)}=S_2^{(n)}(x_2^{(n)})	
\end{matrix}\right.
\end{equation}

\begin{figure}[!hptb]
	\centering
	\includegraphics[scale=.56]{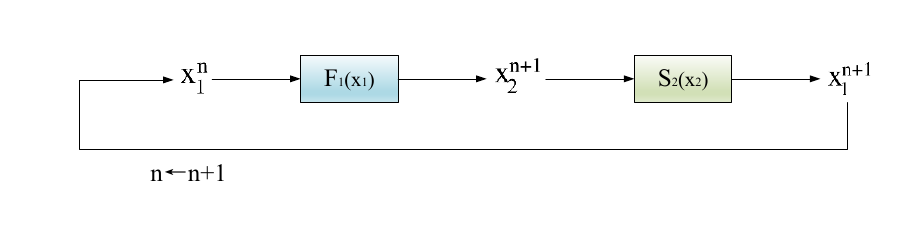}
	\caption{Serial-explicit coupling scheme.}
	\label{fig3}
\end{figure}

Therefore, the explicit coupling scheme is realized through the staggered grid, which allows the solver to have time steps of different sizes, so there are subcycles in the iteration process. preCICE provides special functions to implement these schemes, such as fixed-point iteration and unified adjustment of time steps. The flow of the explicit coupling scheme is shown in Fig.\ref{fig3}.

The principle of explicit coupling is simple and easy to implement, but it also has the disadvantage that it is difficult to overcome, Brummelen\cite{2009Added} describing in detail the features of numerical instability that can not be eliminated even by adjusting the length of time step. Thus, several coupling iterations are performed at each time step until both sides of the solution converged.

\subsubsection{Implicit coupling schemes}
\begin{figure}[!tb]
	\centering
	\includegraphics[scale=.47]{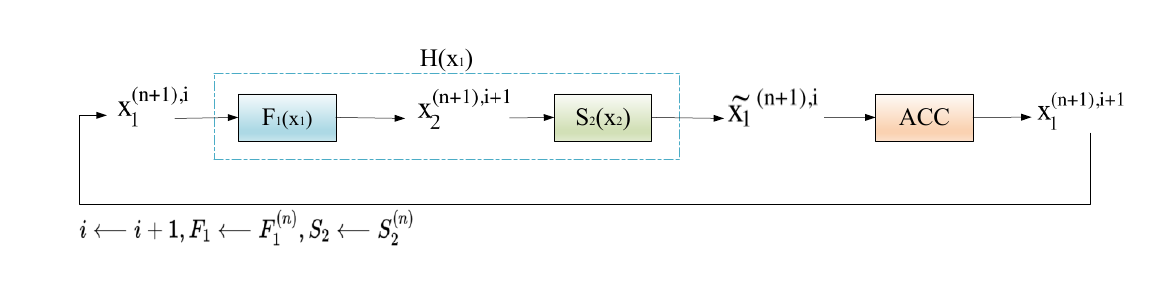}
	\caption{Serial-implicit coupling scheme.}
	\label{fig4}
\end{figure}

The implicit coupling scheme is implemented through a fixed-point iteration, which can ensure the efficiency and precision of the calculations. The form of the fixed point iteration of equation (\ref{ns_eq9}) and (\ref{ns_eq10}) for the serial-explicit coupling scheme can be expressed as equation (\ref{ns_eq12}). The serial-implicit coupling scheme is shown in Fig.\ref{fig4}. 

Where H($x_1$) denotes the solver $F_1$ and $S_2$ are solved couplingly ($H=S_2\odot F_1$), ACC denotes subsequent processing steps, such as acceleration methods.

\begin{equation}\label{ns_eq12}
\left\{\begin{matrix}
x_2^{(n+1),i+1}=F_1^{(n)}(x_1^{(n+1),i})  \\
\\
x_1^{(n+1),i+1}=S_2^{(n)}(x_2^{(n+1),i+1})	
\end{matrix}\right.
\end{equation}

Equation (\ref{ns_eq12}) indicates that both solvers use the value of time $n+1$ when calculating. Similarly, the fixed-point iteration form of equation (\ref{ns_eq11}) for the parallel-explicit coupling scheme can be written as equation (\ref{ns_eq13}).
\begin{equation}\label{ns_eq13}
\left\{\begin{matrix}
x_2^{(n+1),i+1}=F_1^{(n)}(x_1^{(n+1),i})  \\
\\
x_1^{(n+1),i+1}=S_2^{(n)}(x_2^{(n+1),i})	
\end{matrix}\right.
\end{equation}

In each non-convergent sub iteration, preCICE requires that checkpoints are provided and the latest state is saved. If the minimum residual required for convergence is not reached, the sub-loop will reload the latest state and enter the next iteration until the residual requirements are met.
Therefore, in the process of convergence, it still satisfies the serial-coupling scheme.
\begin{equation}\label{eq_pre01}
x_1^{(n+1)}=S_2^{(n)}(x_2^{(n+1)})=S_2^{(n)}(F_1^{(n)}(x_1^{(n+1)}))
\end{equation}    

Therefore, the implicit coupling scheme is the process of solving the fixed-point equation (\ref{eq_pre02}) for the type of serial or Gauss-Seidel.
\begin{equation}\label{eq_pre02}
x_1^{(n+1)}=S_2^{(n)}\odot F_1^{(n)}(x_1^{(n+1)})
\end{equation}  

For the parallel implicit-coupling scheme, we have
\begin{equation}\label{eq_pre03}
\begin{pmatrix}
x_1^{(n+1)} \\
\\
x_2^{(n+1)}
\end{pmatrix}=\begin{pmatrix}
0  & S_2^{(n)}\\
\\
F_1^{(n)}  & 0
\end{pmatrix}\begin{pmatrix}
x_1^{(n+1)} \\
\\
x_2^{(n+1)}
\end{pmatrix}
\end{equation} 

In summary, equations (\ref{eq_pre02}) and (\ref{eq_pre03}) are fixed-point equations that are decoupled by an implicit-coupling scheme. Meanwhile, preCICE provides a variety of methods to solve them, such as underrelaxation, adaptive Aitken underrelaxation and various sophisticated Quasi-Newton solver methods. All these methods can be used to solve any type of fixed-point equations, and more detailed solving steps are described in literature\cite{bungartz2016precice}.

\subsection{Data mapping}
\begin{figure*}[!htb]
	\centering
	\includegraphics[scale=.8]{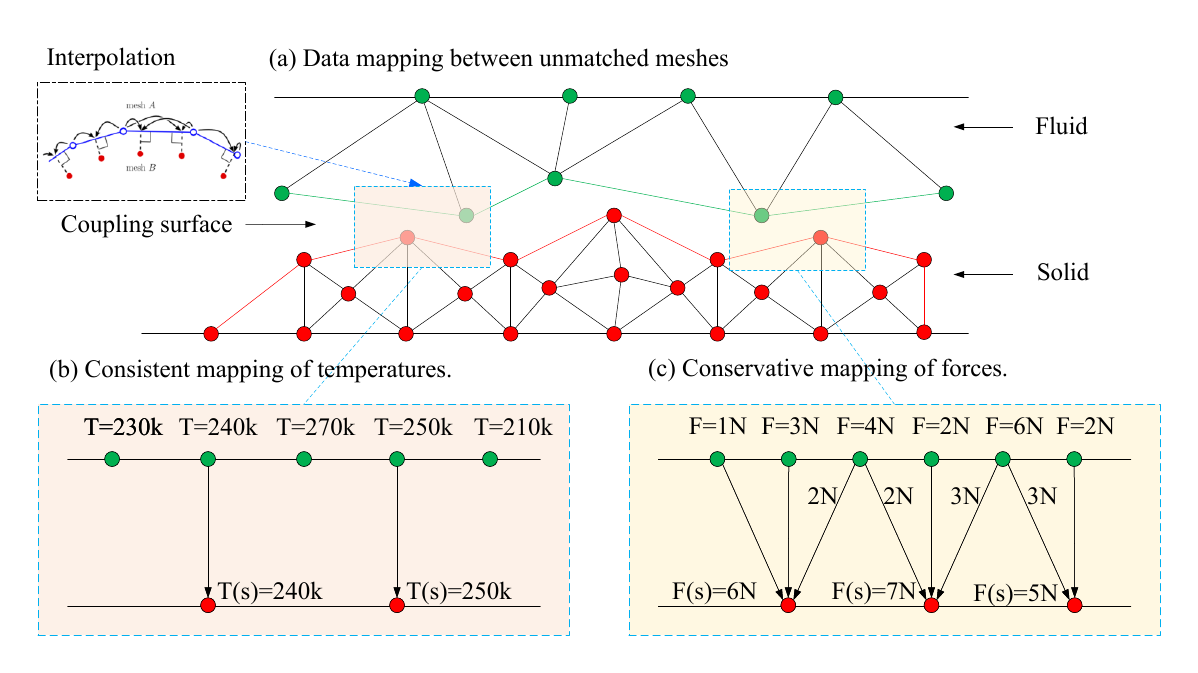}
	\caption{(a) represents the mismatched mesh coupling surface of fluid and solid, (b) represents the consistent mapping principle of temperature $T$, and (c) describes the conservative mapping principle of force $F$.}
	\label{fig5}
\end{figure*}

preCICE supports the coupling of multiple solvers. Generally speaking, different solvers adopt different meshes on the coupling surface (many-to-many relationships). It requires data mapping of unmatched meshes on the coupling surface, and this mapping relationship can not destroy the basic conservation equations such as mass conservation and energy conservation, otherwise, the simulation results are difficult to converge. Therefore, preCICE provides consistency and conservatism principles to deal with data mapping between unmatched grids of Fig.\ref{fig5} (a).

Consistent mapping and conservative mapping are used in different situations. First of all, the consistency form just copies the data from one grid to another grid of the receiver, and its value does not change. For example, in heat transfer simulation, the change of temperature in a certain region is consistent. As shown in Fig.\ref{fig5} (b), the temperature between adjacent grids should use the consistent mapping principle, which ensures that the corresponding grid temperature within the same area is consistent (the grid is the basic unit of atomic operation).
In addition, Newton's third law states that the force and its reaction are equal in magnitude. For solving Fluid-Structure Interaction problem, the transfer of force should strictly follow the conservative mapping principle. As shown in Fig.\ref{fig5} (c), the resultant force on both sides of the fluid and solid should remain unchanged. Therefore, the forces between adjacent mesh interpolation points shall be evenly distributed, such as $F=4N=2N+2N$ and $F=6N=3N+3N$.

preCICE provides the following methods and their variants to implement the data mapping function, which can ensure the conservation conditions between the coupling surfaces.
\begin{itemize}
	\item \textbf{Nearest neighbor:} Its principle is to find the nearest mesh points for interpolation in space. It does not need any topology information. This method is simple and easy to implement, but its disadvantage is that it only has first-order accuracy.
	\item \textbf{Nearest projection:} In this method, the target grid points are projected onto the source grid elements according to the grid topology information, and the source grid elements are linearly interpolated according to the coordinates of the projection points. Because the orthogonal distance between the coupling interface and the grids on both sides is very small, therefore, this method has second-order accuracy.
	\item \textbf{Radial Basis Function (RBF):} This method does not need either grid topology information or projection information. It constructs an interpolation function (RBF) centered on the source grid points. Through this function, the target grid points can be interpolated and the interpolation results can be obtained. Although RBF is globally supported, considering the complexity and efficiency of calculation, its application in practice is only limited to a small range. At the same time, preCICE also provides a variety of interpolation functions that can be consulted in literature\cite{bungartz2016precice}, among which Gaussian and Thin Plate Spline are the most widely used RBF interpolation functions.
\end{itemize}

\subsection{Data communication}
It is a great challenge to couple multiple solvers with parallel computing for efficient communication. Each participant may be executed by multiple processes and distributed on multiple nodes of the supercomputer. Therefore, preCICE provides MPI technology to solve the problem of large-scale parallel communication. It allows open-source programs to run efficiently across nodes. However, for closed-source commercial software, no corresponding interface is provided for calling. preCICE provides a lower-layer TCP/IP sockets communication protocol to support this function which is difficult to program. Each participant communicates point-to-point in parallel, and then it sets up a master process to steer other slave processes. These communication relationships are set at the time of initialization. They are configured through an XML file. preCICE will automatically recognize the mesh decomposition and establish contact with the corresponding processes. This communication method also has some disadvantages, for example, it is difficult to deal with adaptive mesh and FSI problems of fluid immersed boundaries.

\subsection{List of coupled solver}
When using preCICE for coupling development, each participant (solver) needs to make appropriate modifications. Among them, the solver interacts with preCICE through the adapter, as shown in Fig.\ref{fig2}. Therefore, the corresponding adapter for fluid or solid must be developed to use preCICE for coupling. Although preCICE was written in C++. However, the interface it provides also supports C, Python, Fortran and other widely used programming languages. This highly modular software architecture enables it to be flexibly compatible with most of the current influential CFD/CSD software, including open-source and closed-source commercial software. At present, the compatible open source software has been released on the official website (\href{https://precice.org/adapters-overview}{www.precice.org/adapters-overview}).

At present, preCICE has integrated many open/closed-source software including those enumerated in Table \ref{tab01}, and it has achieved many successful applications in FSI \cite{rodenberg2021fenics}, conjugate heat transfer \cite{kim2020simulation}, multi-phase flow \cite{bungartz2011}, and acoustic fields \cite{bungartz2016}. It has a strong function and rich development documentation, at the same time, software community users very active support the preCICE updating constantly, and it attracts more and more researchers worldwide to improve and strengthen the function of preCICE. In future work, the development team will provide the support to change the interface grid.

\section{Particle-mesh coupling method}\label{section4}
In this section, we will give a comprehensive introduction to the particle-mesh coupling (PMC) method proposed in this paper. In the section \ref{section3}, the theoretical basis and application of preCICE are described in detail. It provides strong coupling capabilities, making it the preferred choice for coupled development. However, it only provides method with mesh-based coupling (eg. unstructured or ALE meshes), up to now, no adapter has been released on the official website to support meshless coupling. In addition, the meshless method has inherent advantages in dealing with large deformation, high-speed impact response of materials, and free surface flow. It is very suitable for use as a solver for simulating fluids in preCICE. Therefore, developing a meshless adaptation method that is compatible with preCICE is an important topic.

This paper presents an adaptation method based on SPH and preCICE for simulating FSI problems, which supports the direct interpolation between particle nodes and mesh elements, and then it couples the particle (or meshless) method to the preCICE surface coupling framework.
\begin{figure*}[!htb]
	\centering
	\includegraphics[scale=.4]{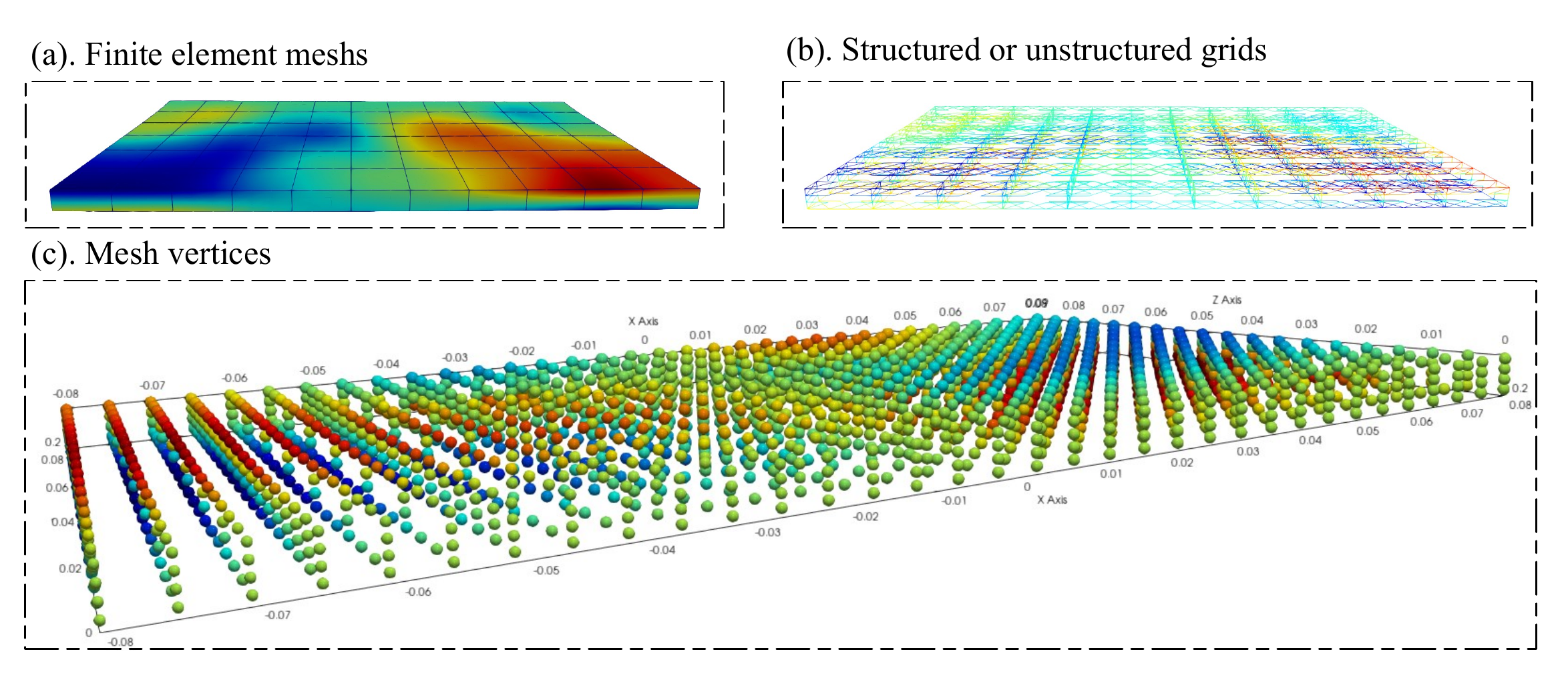}
	\caption{(a) represents finite element mesh, and (b) represents structured/unstructured mesh. (c) corresponds to the vertex of mesh (b). Where the color indicates that there are changes in physical quantities (eg. temperature, stress and displacements) in the mesh/node.}
	\label{fig6}
\end{figure*}

In traditional mesh methods, for example, FEM is often used to discretize solid structures. Fig.\ref{fig6} (b) shows the structural/unstructured mesh after discretization of the solid structure, and Fig.\ref{fig6} (c) shows the vertex of the mesh, which also has a dependency on edges. preCICE couples only according to the mesh vertices on the contact. For example, the physical variables (eg. displacements) of the mesh vertices of the six faces in Fig.\ref{fig6} (c) are usually used as the input of preCICE, and then the preCICE at the other part also gives the output of the corresponding vertices (eg. forces). 
This interaction method becomes very concise for the processing of coupling boundary, which only needs to send the physical information of mesh vertices. Thanks to various powerful data mapping functions provided by preCICE, it does not even need to consider the adjacent relationship between them, and only deals with Dirac boundary conditions (first-order displacements and forces).

\begin{figure*}[!thb]
	\centering
	\includegraphics[scale=.32]{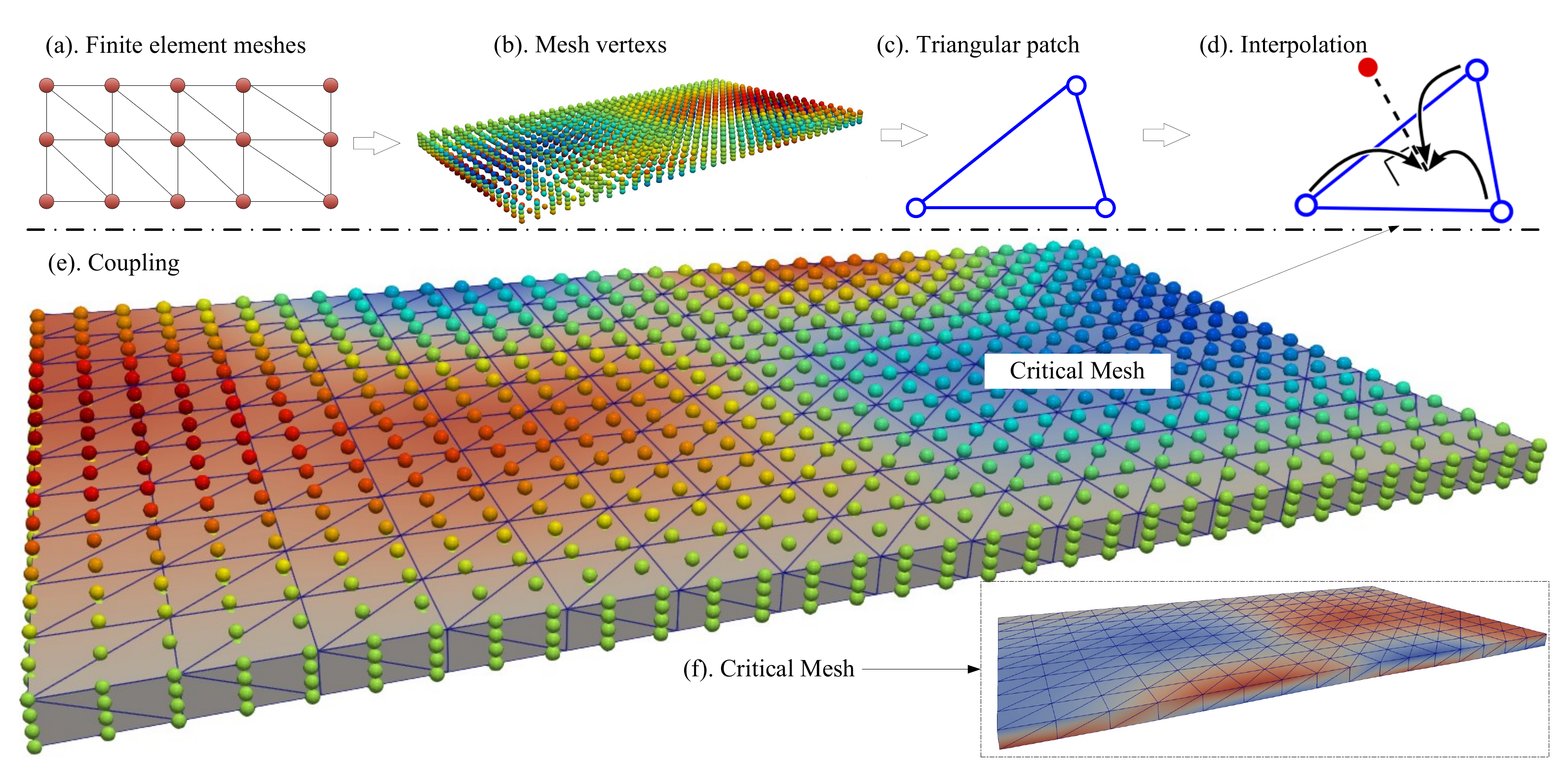}
	\caption{Critical grid provides an interpolation method based on structural/unstructured mesh, which is composed of many triangles/quadrilaterals. The finite element mesh interpolates the critical grid according to the vertex form. Where the color represents the change of physical quantities at the position, such as forces, temperature and displacements.}
	\label{fig7}
\end{figure*}

As a pure Lagrangian particle method, the SPH method has no topological relationship between particle nodes. Unfortunately, preCICE provides interfaces that only support mesh-based methods as shown in Fig.\ref{fig6}. 
and these mesh contact relationships are established at the time of initialization and do not change in subsequent calculations, so particle nodes are not compatible with the definition criteria for data mapping interfaces to perform data exchange. To solve this problem, we have introduced a critical grid as the medium for interpolation between particles and grids, the principle of which is shown in Fig.\ref{fig7}.

As shown in Fig.\ref{fig7}, in order to reuse the interface function provided by preCICE, it needs to artificially add a critical mesh as the intermediate medium for data exchange. This critical grid needs to meet the following characteristics. First, the vertex coordinates of the critical grid are located on the Fluid-Structure Interaction interface, and there is only one layer of grid refer to Fig.\ref{fig7} (f), which is an intermediate buffer for data interpolation as shown in Fig.\ref{fig7} (d). Its definition completely follows the standard form of mesh method (vertex coordinates and edges), and it can be regarded as a custom mesh boundary, which is composed of boundary IDs. Secondly, the data structure of grid boundaries and the dependency of edges are always completely unified with the construction method provided by preCICE, and the purpose is to reuse the mesh-based interface provided by preCICE. Finally, preCICE is only the coupling of the contact surface. As the critical area of data exchange between fluid and structure, it belongs to the wet interface. In essence, the critical grid is also a special contact surface. Taking the finite element mesh as an example, the vertices of the finite element meshes and the critical grids are interpolated as shown in Fig.\ref{fig7} (e), both of which are coupled based on the mesh method (the critical grid has only one layer of grid). Therefore, by introducing the critical mesh, the fluid particles that are located in the wet interface can initially call the interface function of preCICE.

In order to solve the coupling problem of the particle (meshless) method, a more flexible coupling framework is established by introducing the critical grid method, so that the solver in the wet interface have more choices. Both the mesh-based method and the pure Lagrangian particle method can provide strong compatibility. preCICE interacts with the data between the boundary grid and the critical grid, and then the coupling process successfully establishes the connection between the fluid and the structure by jointly calling the interface function provided by preCICE, this procedure is completed by the initialization function of preCICE.

\begin{figure}[!thb]
	\centering
	\includegraphics[scale=.2]{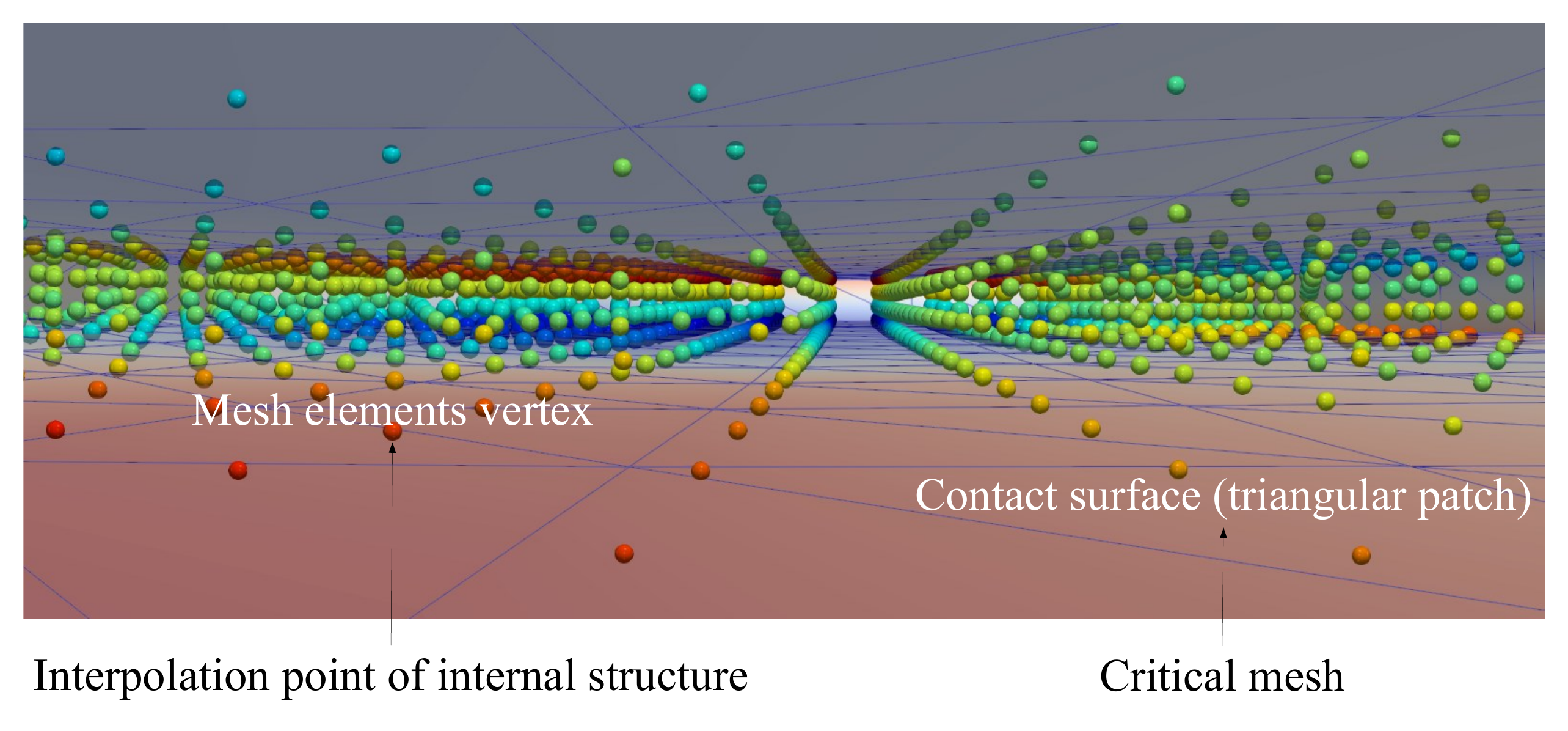}
	\caption{The solid structure part is divided into internal region and external region. The small sphere here represents the mesh vertex of the finite element, and their outside is surrounded by the critical mesh. Among them, the color represents the information corresponding to the change of physical quantities.}
	\label{fig8}
\end{figure}

After introducing the critical grid, as shown in Fig.\ref{fig8}, the solid structure is divided into the internal region and external region. The internal region is calculated by the inherent discrete rules of the FEM method, such as solving elastic and inelastic equations. For the boundary part (external region) of the solid structure, it is not only affected by its internal solution method, but also needs to map the data with the critical grid. The mesh boundary of the finite element will be interpolated with the critical mesh composed of triangular (or rectangular) patches. And the interpolation also needs to meet the law of physical conservation, otherwise, it is easy to cause numerical divergence and calculation errors. Therefore, the solution of solid boundary will be subject to more constraints.

\begin{figure*}[!htb]
	\centering
	\includegraphics[scale=.48]{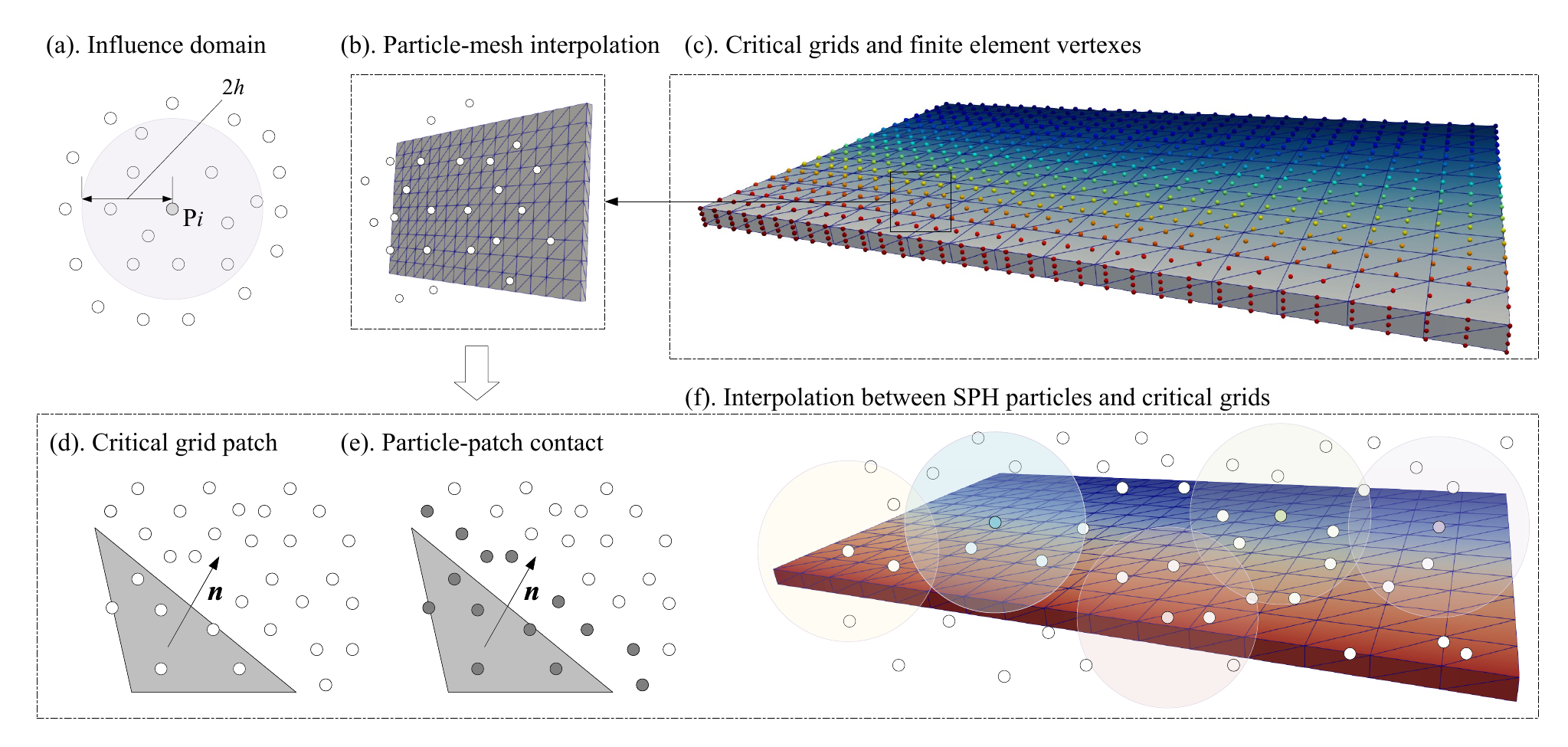}
	\caption{In addition to the data mapping between the critical mesh and the finite element through preCICE, it also needs to interpolate with the particles of SPH. Where (a) represents the influence region of SPH particles, (e) represents the contact relationship between triangular/rectangular patches of critical mesh and particles, $\bm{n}$ is the normal vector of triangular/rectangular patches, and (f) represents the interpolation process between particles and critical mesh.}
	\label{fig9}
\end{figure*}

The critical mesh transfers the corresponding physical quantities to the mesh nodes on the finite element surface through the mapping relationship provided by preCICE. For example, in the FSI problem, the critical grid obtains the force of fluid particles and applies it to the solid structure. preCICE establishes a mapping relationship between the critical grid and the solid surface boundary grid with the form of vertex coordinates. As an intermediate medium, the critical grid applies the force at the fluid end to the solid structure and successfully causes the deformation of the solid structure. The deformation of the solid is transferred to the critical grid through preCICE, and then the deformation of the solid structure is fed back to the fluid particles through the critical grid. It can be seen that by completing multiple interpolations with the critical grid, the fluid particles can achieve indirect contact with the solid grid. The FSI problem based on the interpolation principle of the critical grid is shown in Fig.\ref{fig9}.

Previously, we have designed the docking method of the critical grid and solid grid by calling the external interface function provided by preCICE. However, the particle method can not interact directly with preCICE, so we use SPH particles and critical mesh interpolation to couple the particle method to preCICE.
We call this interpolation process that supports coupled meshes as the particle-mesh coupling (PMC) method. Since it is not restricted by the background mesh (the critical mesh is generated based on the surface of the solid mesh), the theoretical principle of this method allows it to support almost all particle methods to couple with preCICE.

In solving the FSI problem, we use the particle (meshless) method to discretize the fluid part. the data interaction between the critical grid and finite element can be completed through preCICE as shown in Fig.\ref{fig9} (c), in which fluid particles need to interpolate and exchange data with the critical grid (mesh). It requires the force generated by particles to act on the critical mesh. On the one hand, because the fluid particles of the SPH method have an influence region with a radius of 2$h$, the particles of the patch in the critical grid whose distance from the center in the normal vector direction is less than 2$h$ are classified as contact particles. In Fig.\ref{fig9} (e), the contact particles are represented by gray circles, which need to be interpolated with the critical grid, so as to transfer the force generated by the fluid to the solid structure through the critical grid. It should be noted that the particles are considered to have no contact relationship when the distance between particles and patches is greater than 2$h$, and the particles in this part do not need to exchange data with the critical grid. On the other hand, the solid structure takes the received force as the input and calculates the corresponding deformation, and feeds back the deformation results to the critical grid. Then, the critical grid interpolates the displacement with the contacting fluid particles, and finally updates the position information of the fluid particles according to the interpolation results. And their coupling procedure for Fluid-Structure Interaction (FSI) problems can be expressed as follows.
\begin{equation}\label{cm_eq01}
stage1:F_{particles}\overset{forces}{\rightarrow}  \Upsilon_{medium}  \overset{forces}{\rightarrow} \Psi_{coupling}  \overset{forces}{\rightarrow}\ S_{meshes}
\end{equation}

\begin{multline}\label{cm_eq02}
stage2:F_{particles}\overset{displacements}{\leftarrow}  \Upsilon_{medium}  \overset{displacements}{\leftarrow}  \Psi_{coupling} \\
\overset{displacements}{\leftarrow}\ S_{meshes}
\end{multline}

Where $F$ and $S$ represent fluid and solid modules for FSI problems respectively, $\Upsilon$ represents critical grid, and $\Psi$ represents preCICE coupling procedure. To sum up, this is a bidirectional coupling procedure, which includes two stages: the first stage is to transmit the forces, and the second stage is to transmit the displacements. The critical mesh acts as a medium to realize the docking of particle method and mesh method under the preCICE coupling framework. Due to the pure Lagrangian nature of the SPH method, there is no topological relationship between particle nodes. Theoretically, the PMC method proposed in this paper can also be applied to the coupling of most meshless methods.

\section{Implementation}\label{section5}
In this section, we will introduce the specific implementation of the PMC method in solving the FSI problem. The realization of the coupling process is mainly divided into two aspects: 1) calculating the force generated by the fluid particles on the critical grid and 2) processing the displacement from the solid structure, which is transmitted to the fluid particles through the critical grid and updating the position information of the fluid particles. In addition, we also use multithreading technology and the identification of coupled boundary particles to further optimize the computational efficiency of the program.
\subsection{Adapter for SPH}\label{section51}
In section \ref{section3} of this paper, we mentioned that the preCICE coupling library and the solver are connected through adapters. This highly modular approach makes preCICE compatible with most solvers. Similarly, for the SPH method based on pure Lagrange, we use PMC method to develop an adapter connecting preCICE, which is called SPH-Adapter. It uses a critical grid as the medium to interact with preCICE and solid solver, and its process can be simplified as shown in Algorithm \ref{Algorithm00}.
\begin{algorithm}[ht]
	\caption{The implementation architecture of the SPH adapter is based on the PMC method for solving FSI problem.}
	\LinesNumbered 
	\label{Algorithm00}
	\KwIn{critical meshes, configuration parameters of solid materials (Poisson's ratio, elastic modulus, density), time step dt}
	\KwOut{Fluid–structure interaction simulation results}
	precice::SolverInterface precice(``FluidSolver",rank,size)\;
	precice.configure(``precice-config.xml")\;
	precice.setMeshVertices();//critical meshes (grids)\;
	dt $\gets$ precice.initialize()\;
	$//$ \textit{main time loop}\;
	\While{precice.isCouplingOngoing()}{
		precice.readBlockVectorData(displacements)\;
		sph$\_$dt $\gets$ beginTimeStep()\;
		dt $\gets$ min(sph$\_$dt,dt); $//$ \textit{coordinate time steps}\;
		solveTimeStep(dt);   $//$ \textit{SPH fluid solution module}\;
		precice.advance(dt); $//$ \textit{steer coupling}\;	
		precice.writeBlockVectorData(forces)\;
		endTimeStep(); 		 $//$ \textit{write results, increment time}\;
	
	}
    $//$ \textit{coupling ends, releasing resources}\;
    precice.finalize()\;
\end{algorithm}


Lines 1-2 of Algorithm \ref{Algorithm00} are the configuration of preCICE, where \textit{precice-config.xml} is a configuration file with XML format, which contains a lot of configuration information. Including the configuration of fluid and solid, the adopted data mapping method, data communication mode, setting time window and convergence conditions, etc. The configuration file needs to be loaded into a path recognized by the SPH solver. Then, set a critical mesh in line 3 as the medium of particle interpolation, which meets all the standards of the preCICE interface and is composed of only one layer of mesh as shown in Fig.\ref{fig7} (f). Its mesh vertices are loaded into the array in the form of coordinates according to the topological relationship, and preCICE provides the functions of mesh connection and data initialization to support these operations. Meanwhile, the adjacency relationship between grids has been determined at the beginning and will not be changed later. After the configuration is completed, line 4 calls the initialization function, and then preCICE will establish communication between the fluid and solid solver.
SPH is solved in line 10, which is responsible for the calculation of fluid in the FSI problem. Line 7 is to obtain the displacement returned by the solid solver based on the critical grid, and line 12 is to transfer the force generated by fluid particles to the solid structure based on the critical grid. Data transmission and coupling communication only occur in the advance function in line 11. In addition, in line 9, it is necessary to ensure that the time step used between fluid and solid is consistent, which avoids the calculation error caused by the disordered time advance.

In summary, we established a link between the preCICE framework and SPH particles through a critical grid, and we will describe in detail the process by which the interpolation is achieved for the critical grid and particles, the interpolation results and solid structures are achieved for the exchange of force and displacement through preCICE.
\subsection{Fluid force interpolation}\label{section52}
\begin{figure}[!htb]
	\centering
	\includegraphics[scale=.8]{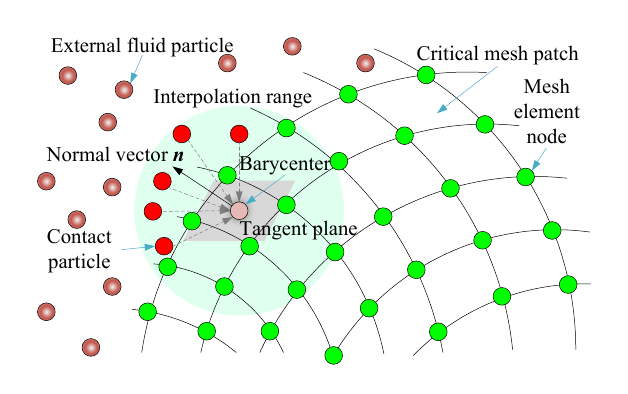}
	\caption{Fluid force (or displacement) interpolation procedure of between fluid particles and critical mesh.}
	\label{fig91}
\end{figure}
In solving the FSI problem, fluid particles transfer the force to preCICE through interpolation with the critical mesh $\Upsilon$, which ensures that the force generated by fluid particles will act on the patch closest to the critical mesh. Therefore, the interpolation between fluid particles and the critical mesh has the following relationship.
\begin{equation}\label{cm_eq03}
\wp_k=\sum_{i=0}^{N-1} \Gamma_{\phi_i}^{\alpha} =\Gamma_{\phi_0}^{\alpha}+\Gamma_{\phi_1}^{\alpha}+\Gamma_{\phi_2}^{\alpha}+\cdots+ \Gamma_{\phi_{(n-1)}}^{\alpha}
\end{equation}

Where $\phi_i$ represents the fluid particle $i$, and $\wp_k$ represents the pressure exerted by the fluid particle on the critical mesh patch $k$. $\Gamma$ represents the force generated by fluid particles in the $\alpha$ dimension ($\alpha=1, 2, 3$). $N$ represents the total number of fluid particles in contact with the patch as shown in Fig.\ref{fig91}.
Obviously, equation (\ref{cm_eq03}) shows that the force on the critical mesh can be obtained by summing all the contact particles, and then the resultant force obtained by particle interpolation is transferred to the solid structure through preCICE as the pressure generated by the fluid on the contact surface. In addition, the particle interpolation method can calculate the resultant force by using two methods. The first method is to calculate directly through the particle node using Newton's second law. First, the acceleration of the particle is obtained through the momentum equation (\ref{basic_eq13}), and then the acceleration and particle mass are used to calculate the corresponding force, Finally, the forces of all particles are added and transferred to the critical grid as the force generated by the fluid on the coupling boundary. However, because this method relies on contact particles as the calculation carrier, there is an obvious disadvantage of this method, which is vulnerable to numerical errors caused by particle distribution disturbances.
In order to make up for this defect, equation (\ref{cm_eq04}) describes the solution method based on particle pressure integration.
It uses the basic principle of integration to solve the contact force generated by fluid particles on solid structures. This integral representation method based on continuous function has the advantages of robustness and stability, and it can better maintain numerical stability in the calculation process.
\begin{equation}\label{cm_eq04}
\wp=\int_{\Omega }P\mathrm{d}S 
\end{equation}
Where $P$ is the pressure calculated from equation of state (\ref{ns_eq7}), $\Omega$ is the range of contact particles, and $\mathrm{d}S$ is the area element. Therefore, equation (\ref{cm_eq04}) shows that the force on the patch $\wp$ can be obtained by solving the integral of the pressure on the contact surface by its contact particles, which is less affected by the particle distribution. Therefore, this calculation method is more accurate than the method based on a single particle to solve the resultant force for critical grid. In order to realize the interpolation of force, we first define the particle contact method (PCM) as shown in algorithm \ref{Algorithm01}.
\begin{algorithm}[ht]
	\caption{The particle contact method (PCM), which realizes the topological connection between particles and critical mesh, and the subsequent force interpolation is implemented based on this correspondence.}
	\LinesNumbered 
	\label{Algorithm01}
	\KwIn{array of all fluid particle $\phi$, critical mesh patch array $\Re$, the threshold of contact distance is $\Theta$, and smooth length $h$.}
	\KwOut{minimum distance set $\chi$ of contact particles and patches.}
	Initialization parameters\;
	$\Theta \longleftarrow 2h$\;
	$//$ \textit{this loop uses multithreaded parallelism}\;
	\For{j in $\phi_{n-1}$}{
		$distMin^j\longleftarrow getDistParticleToPatch(\phi_j,\Re_{start})$\;
		\For{$\Re_{start}<i<\Re_{end}$}{
			$dist_{i}^j\longleftarrow getDistParticleToPatch(\phi_j,\wp_i)$\;
			\If{$dist_{i}^j < \Theta$}{
				\eIf{$dist_{i}^j < distMin^j$}{
					$\chi_j\longleftarrow dist_{i}^j$\;
					$distMin^j\longleftarrow dist_{i}^j$\;					
				}{
					$\chi_j \longleftarrow distMin^j$\;
				}
			}
		}
	}
\end{algorithm}

Algorithm \ref{Algorithm01} describes the contact algorithm between particles and patches, which is used in the interpolation process of particles and critical meshes, including force and displacement interpolation. Among them, the line 2 sets the contact threshold $\Theta$, which is generally set to 2$h$. The fourth line represents the process of searching the nearest critical mesh patch for particles close to the coupling surface. Due to the weak dependence, the module can adopt multi-threaded parallel technology to optimize the efficiency of calculation. In line 7, calculate the distance between particles and patches, which is compared with the contact threshold in line 8. When the distance is less than the set threshold, we judge that there is contact between fluid and solid. On the contrary, when the contact distance is greater than the threshold $\Theta$, it is considered that there is no contact relationship. Finally, in line 10 or 13, we save the particles in contact around the solid to array $\chi$, and it is called contact particles, which are provided for subsequent interpolation with the critical mesh. The particle mesh contact process is described in Fig.\ref{fig9} (d) and (e). It should be noted that the contact particles calculate the shortest distance, which ensures that each particle almost establishes a contact relationship with only one patch. This correspondence is a one-to-many mapping relationship, that is, a patch contains multiple contact particles, but a contact particle only corresponds to the target patch. Therefore, based on this contact relationship, the interpolation of particle force can be implemented by the following algorithm \ref{Algorithm02}.
\begin{algorithm}[ht]
	\caption{The particle and the critical mesh are interpolated by force, which transfers the force generated by the particle on the coupling surface to the solid structure through preCICE in the standard form of the critical mesh.}
	\LinesNumbered 
	\label{Algorithm02}
	\KwIn{critical mesh patch array $\Re$, contact particle set $\chi$, the acceleration of fluid particle is $ace_{\phi}$, the mass of fluid particles is $mass_{\phi}$.}
	\KwOut{resultant force $\wp$ on critical mesh patch.}
	Initialization parameters\;
	$//$ \textit{this loop uses multithreaded parallelism}\;
	\For{$\Re_{start}<i<\Re_{end}$}{
		$\wp_i\longleftarrow 0$\;
		\For{$j$ in $\chi$}{
			update equation (\ref{cm_eq03})\; 
			$ \Gamma_j^i\longleftarrow mass_{\phi_j}\cdot ace_{\phi_j}$ or update equation (\ref{cm_eq04})\;			
			$\wp_i\longleftarrow  \wp_i+\Gamma_j^i$\;
		}
	}
\end{algorithm}

Algorithm \ref{Algorithm02} describes the simplified process of particle force interpolation, in which line 3 traverses all patches $\Re$ of the critical mesh. Then, according to the contact relation $\chi$ provided by algorithm \ref{Algorithm01}, we calculate the acceleration of particles according to the momentum equation, then calculate the force generated by a single particle on the critical grid according to Newton's second law (the particle mass is a constant), and finally calculate the resultant force $\wp$ generated by all particles on the critical grid as a whole. Then the resultant force $\wp$ is transferred to the solid structure through preCICE. At the same time, the force interpolation can also be implemented through equation (\ref{cm_eq04}).  In theory, although it has high accuracy, the disadvantage is that it will be difficult to realize the pressure interpolation on the complex coupling surface, and this process will consume additional computing resources. To sum up, after realizing the critical mesh and particle contact algorithm, the force of particles on the solid mesh can be calculated, and the accuracy of interpolation will also have a certain impact on the accuracy of the final coupling results. Therefore, it is very important to choose an efficient particle mesh interpolation algorithm to obtain high-precision results.

\subsection{Displacement interpolation}\label{section53}
In section \ref{section52}, we describe the force interpolation procedure of particles, and then the solid structure calculates the displacement of deformation according to the force exerted by particles on its surface, and the displacement is transmitted back to the critical grid through precCICE. Therefore, it is necessary to interpolate the displacement received by the critical mesh to the contact particles and update the position coordinate information of the particles. The displacement interpolation procedure can be described by algorithm \ref{Algorithm03}.
\begin{algorithm}[ht]
	\caption{The particle and the critical mesh are interpolated by force, which transfers the force generated by the particle on the coupling surface to the solid structure through preCICE in the standard form of the critical mesh.}
	\LinesNumbered 
	\label{Algorithm03}
	\KwIn{critical mesh patch array $\Re$, contact particle set $\chi$, spatial dimension $\alpha$, particle velocity $V$, coordinate position $P$ of particles, the absolute displacement of the critical grid relative to the initial position is $\xi$.}
	\KwOut{the particles information of position and velocity  after correcting.}
	Initialization parameters\;
	$//$ \textit{this loop uses multithreaded parallelism}\;
	\For{$\Re_{start}<i<\Re_{end}$}{
		$\xi_{relative}\longleftarrow \xi_{last}-\xi_{current}$\;
		\For{$j$ in $\chi$}{
			$//$ \textit{correct the position information of contact particles}\;
			$ P[\alpha][j]_{correction}\longleftarrow P[\alpha][j]+\xi_{relative}$\;		
			$//$ \textit{prevent particles from non-physical penetration of boundaries for velocity correction using equation (\ref{cm_eq05})}\;	
			$V[\alpha][j]_{correction}\longleftarrow Correction(V[\alpha][j])$\;
		}
	}	
\end{algorithm}

Algorithm \ref{Algorithm03} describes the interpolation procedure of displacement between the critical grid and its corresponding contact particles. It includes two parts: 1) the calculation of critical grid displacement and 2) the correction and update of particle position information. It is necessary to calculate the last displacement $\xi_{last}$ of the critical grid and the current displacement $\xi_{current}$ in line 4, and correct the current contact particles according to the relative displacement $\xi_{relative}$ in line 7 to prevent the boundary particles from non-physical penetration of the coupling surface. However, the contact particles with high speed may still penetrate the coupling boundary in the system update stage of SPH (the information of physical quantities such as particle position and speed will be updated). Therefore, it is also necessary to correct the speed, it needs to satisfy the boundary conditions of equation (\ref{cm_eq05}).
\begin{equation}\label{cm_eq05}
\left ( \frac{\mathrm{d} \bm{\xi} }{\mathrm{d} \bm{t}}-\bm{v}  \right ) \cdot \bm{n}=0
\end{equation}   

Where $\bm{\xi}$ and $\bm{v}$ represent the displacement and velocity of the critical mesh respectively, and $\bm{n}$ represents the normal vector at the critical mesh patch. Equation (\ref{cm_eq05}) indicates that to prevent fluid particles from penetrating the boundary at the wet interface, the velocity of the contact particles in the normal vector direction of the critical mesh patch needs to be less than or equal to the moving velocity of the patch, the velocity of the contact particles is corrected using equation (\ref{cm_eq05}) in line 9. At the same time, in order to improve computing efficiency, the procedure can use multi-threaded acceleration technology in line 4. To sum up, the displacement interpolation needs to deal with the contact relationship between particles and the critical grid, avoid the non-physical penetration of the contact particles to the coupling boundary resulting in calculation errors.

\section{Numerical examples}\label{section6}
In this section, we will give numerical examples to verify the performance of the particle mesh adaptation method for solving the problems of FSI. In this experiment, we use the SPH method to solve fluid, and its solver is released on GitHub\footnote{\href{https://github.com/GabrielDigregorio/SPH\_method}{www.github.com/GabrielDigregorio/SPH\_method}} and can be downloaded freely. In addition, we use the finite element method to solve the solid structure, and its solver adopts Deal.II\footnote{\href{https://www.dealii.org/}{www.dealii.org}}, and the adapter of Deal.II connecting preCICE can be obtained on the official website\footnote{\href{https://github.com/precice/dealii-adapter}{www.github.com/precice/dealii-adapter}}. 
\subsection{Breaking dam flow impact on an elastic plate}
\begin{figure*}[!htb]
	\centering
	\includegraphics[scale=.8]{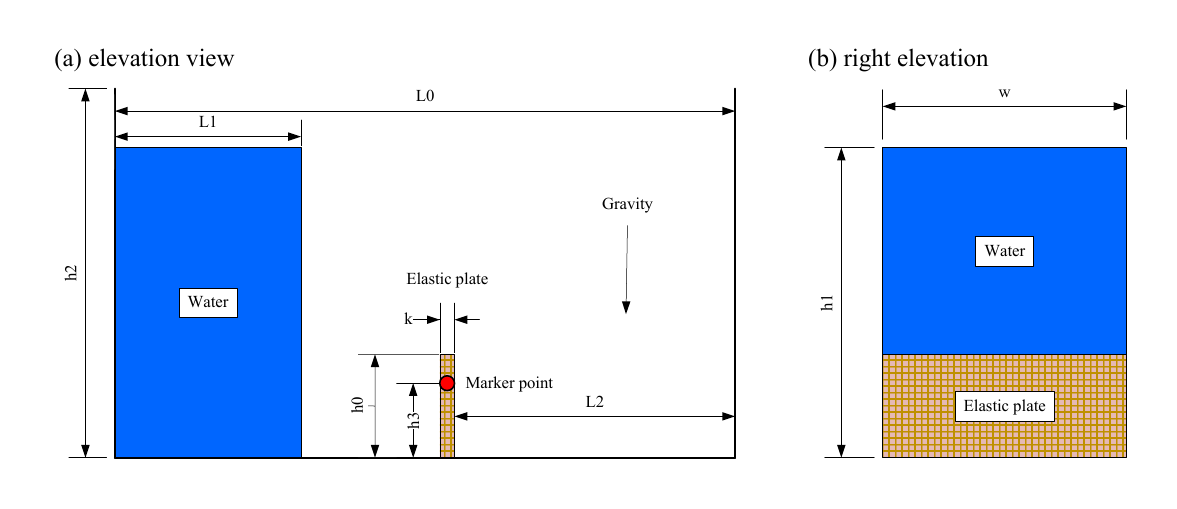}
	\caption{The dimensions of the elastic plate, water, and water tank are set for the dam breaking flow model, where (a) and (b) are the elevation view and the right elevation, respectively.}
	\label{fig12}
\end{figure*}
This example gives the numerical simulation results of the classical FSI problem of the dam breaking flow impacting the elastic plate to verify the accuracy of the PMC method, the benchmark was completed by Liao \textit{et al.}\cite{liao2015free}. And it is an important example to verify the FSI problem. The setting of the model is shown in Fig.\ref{fig12}, in which the bottom of the elastic plate is fixed and made of rubber, and its corresponding parameters are as follows, the density is $\rho=1161.54$ kg/m$^3$, the Young's modulus is $E=3.5\times 10^6$ Pa, and the Poisson's ratio is $\nu=0.45$, the model is in a gravity field with a vertically downward direction, and the gravity acceleration is $\bm{g}=9.8$ m/s$^2$. The size of the model is shown in Table \ref{tab02}.
\begin{table*}[thp]  
	\centering
	\caption{Geometric parameters of dam breaking flow impacting elastic plate (unit: meters).}  
	\label{tab02}
	\begin{tabular}{p{1.5cm}p{1.5cm}p{1.5cm}p{1.5cm}p{1.5cm}p{1.8cm}p{1.8cm}p{1.3cm}p{1.3cm}}
		\toprule
		Length of water tank $L_0$ (m) & Height of water tank $h_2$ (m) & Width of water $w$ (m) & Length of water $L_1$ (m) & Height of water $h_1$ (m)& Height of elastic plate $h_0$ (m)& Thickness of elastic plate $k$ (m)& Margin $L_2$ (m) & Marker point $h_3$ (m)\\
		\midrule
		0.8 & 0.6 & 0.12 & 0.2 & 0.4 & 0.1 & 0.004 & 0.2 & 0.087\\
		\bottomrule
	\end{tabular}
\end{table*}

\begin{figure*}[!tb]
	\centering
	\includegraphics[scale=.4]{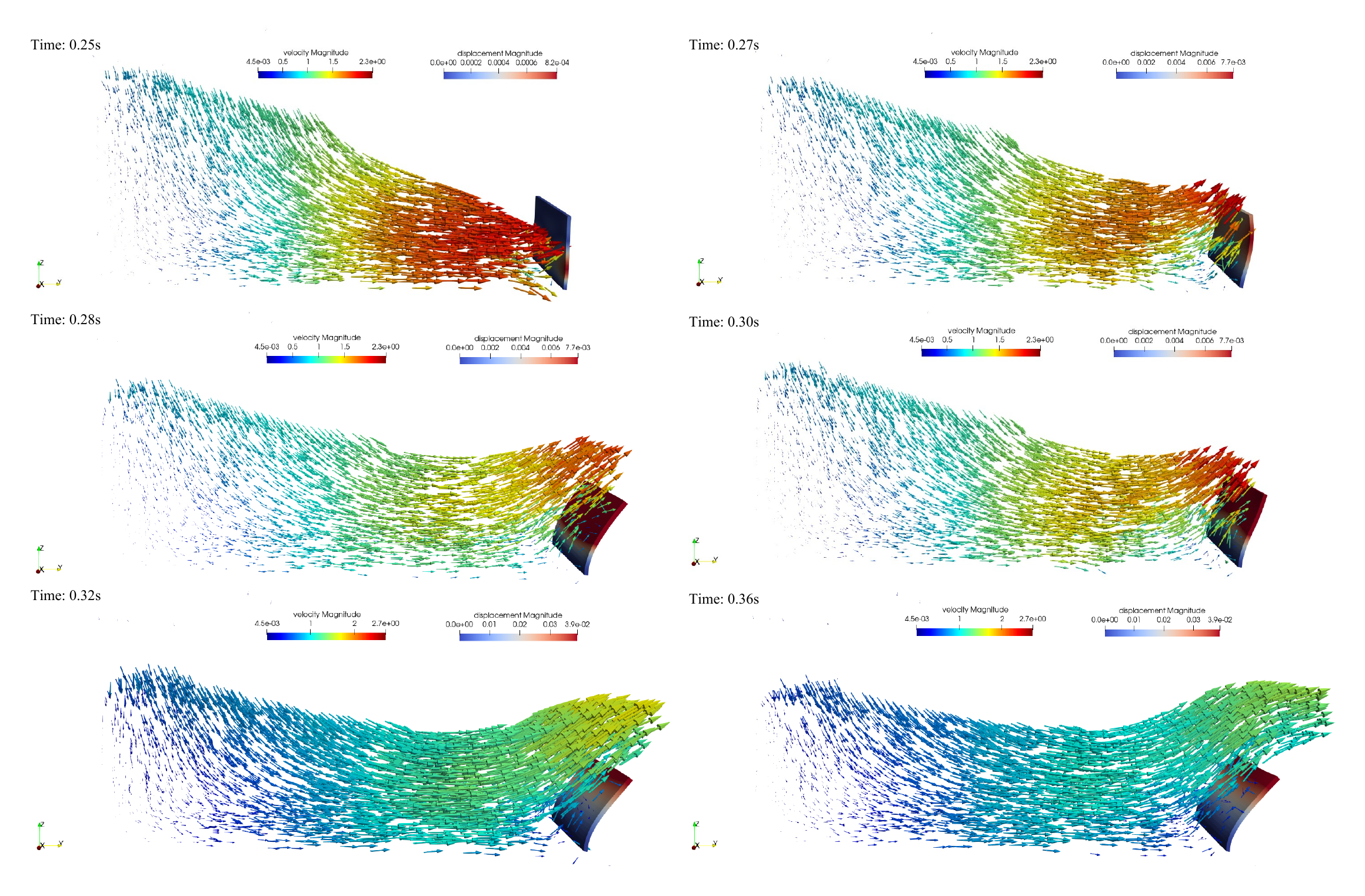}
	\caption{The flow field of the dam break flows at different times and the deformation procedure is caused by the impact on the elastic baffle.}
	\label{fig130}
\end{figure*}

In this experiment, the width of the elastic plate is 0.2 m, the thickness is 0.004 m, and the height is 0.1 m. The spacing between the rubber plate and the right wall is 0.2 m, the length of the water is 0.2 m, and the height is 0.4 m, the rubber plate was marked with points, where $h_3$=0.087 m. The whole device is placed in a water tank with the length-width-height of $0.8\times0.2\times 0.6$. This example is used to verify the accuracy of the numerical results of the PMC method. The dam breaking flow is solved by the SPH method, and the elastic baffle is solved by the FEM method. The simulation results are shown in Fig.\ref{fig130}.

\begin{figure*}[!htb]
	\centering
	\includegraphics[scale=.5]{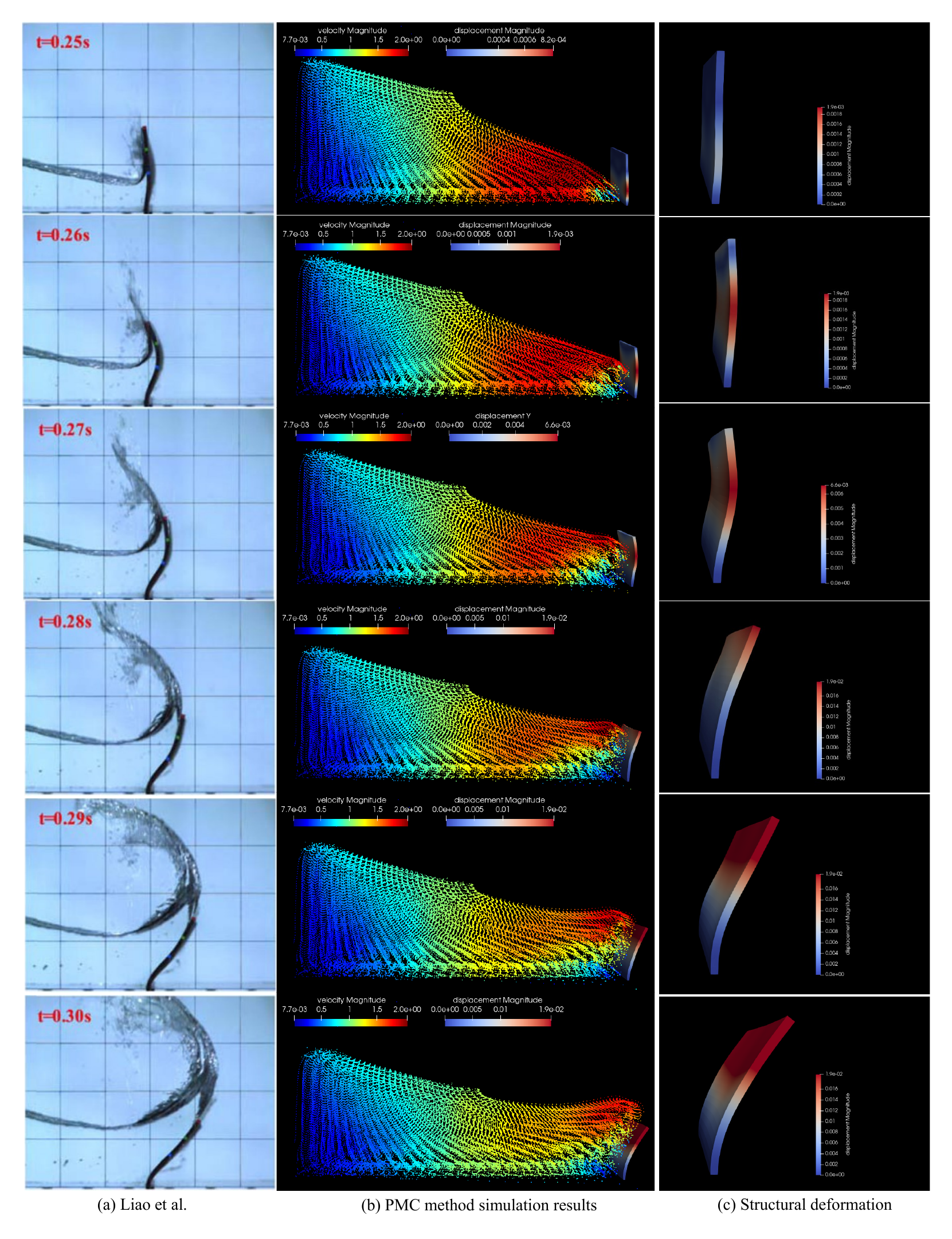}
	\caption{Compare the PMC Method to simulate the dam break flow impacting the elastic plate and the actual experiment of \cite{liao2015free}, the time is from 0.28 s to 0.30 s.}
	\label{fig132}
\end{figure*}
Fig.\ref{fig130} describes the impact process of dam break flow on the baffle, which can better verify the performance of the PMC method. It records the flow field and the deformation process of the baffle within 0.25 to 0.36 seconds. The color of the fluid indicates the flow velocity of the flow field, and the color of the elastic baffle indicates the deformation. First, within 0.25 s, the fluid contacted the elastic plate and began to impact it. When the time reaches 0.27 seconds, the elastic baffle begins to produce concave deformation due to the impact of water flow. Obviously, the fluid particles contacting the elastic baffle do not produce non-physical penetration to the baffle, which is the result of using equation (\ref{cm_eq05}) to constrain the boundary conditions. When the time is within 0.32 s, because the bottom of the baffle is fixed, its shape variable is smaller than that of the top. Therefore, the water velocity near the bottom of the baffle will drop sharply and the direction will also change. At the same time, the top of the baffle is movable, and it will produce greater deformation after being impacted by the water flow. And most of the fluid particles will continue to flow forward across the top of the baffle. In the time from 0.32 to 0.36 s, the flow velocity of the flow field has a large loss until it is finally reduced to 0. Finally, the elastic baffle will return to its original position due to the force balance. The contact radius we used in this simulation example is $\Theta=0.01$ m. Different contact radii will affect the simulation results. For example, if the contact radius is too small, even if the boundary conditions are imposed, the fluid particles can still produce non-physical penetration to the elastic baffle, which is caused by the updating of particles in the calculation process of the SPH method. In this way, the faster particles penetrate into the interior of the solid without modifying the coupling boundary conditions, which is illegal. Therefore, in the actual simulation, we need to comprehensively consider the impact of various factors on the simulation results.
\begin{figure*}[!htb]
	\centering
	\includegraphics[scale=.7]{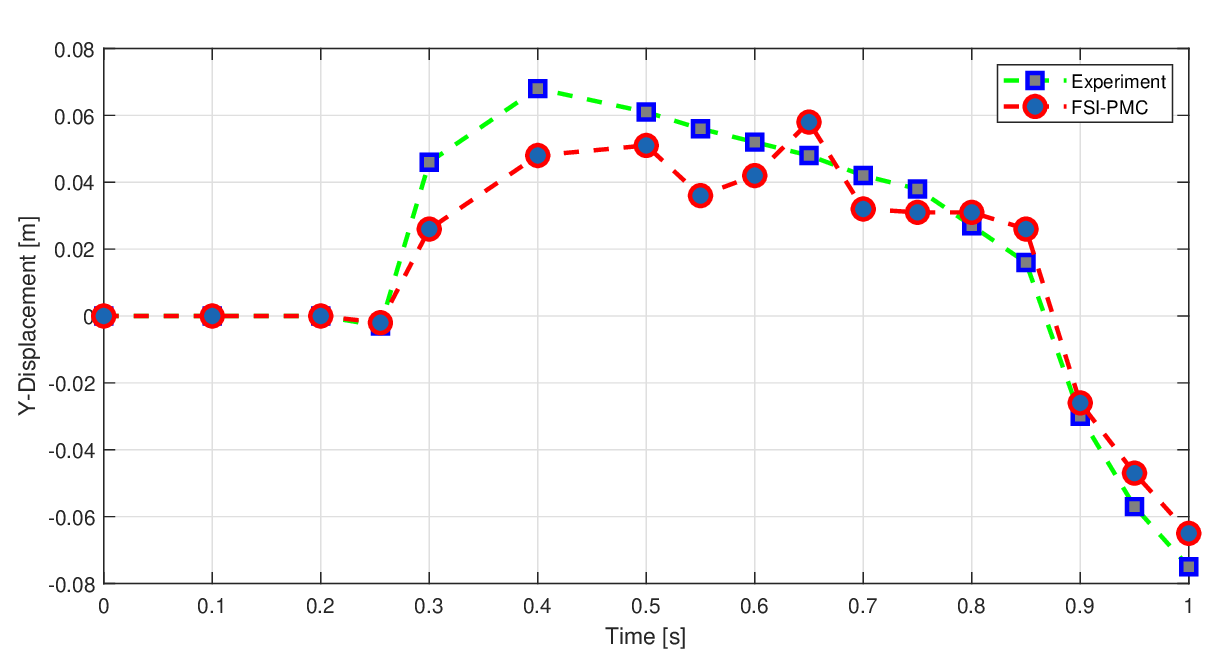}
	\caption{Comparison between the displacement along the Y-direction of the marker point and the actual experiment of Liao \textit{et al.}\cite{liao2015free}, the time is from 0 s to 1 s.}
	\label{fig135}
\end{figure*}

In order to further verify the performance of the PMC method, we compared it with the actual experimental results of Liao\cite{liao2015free}. As shown in Fig.\ref{fig131} and Fig.\ref{fig132}, the experiment gives the results within 0.25 s to 0.3 s, and the flow procedure at different times is recorded by a high-speed camera. The fluid particles first touched the baffle within 0.25 s, and the structural module has been deformed. With the further increase of the fluid velocity, the structure part has depression, and the simulation results are consistent with the actual experimental results. The structure (elastic plate) is further deformed within 0.28 s after being impacted by water flow. Finally, the structure was deformed greatly within 0.3 s, and the fluid particles passed through its top. However, the experimental results show that some fluid particles are thrown out, and there are subtle differences in the location of local fluid. The reason is that the actual experiment also has an influence on airflow. 
In addition, as shown in Fig.\ref{fig135}, the displacement in the Y-direction is marked (refer to the marker point in Fig.\ref{fig12}), it can be seen that this procedure describes the deformation process of the marker point simulated by the PMC method within one second, and it is compared with the actual experimental results of Liao \cite{liao2015free}. Firstly, the fluid particles come into contact with the rubber baffle within 0.25 s and force it to produce deformation displacement along the Y-direction. Secondly, when the simulation time reaches 0.4 s, the deformation displacement generated by the baffle reaches its maximum. Overall, between 0.3 s and 0.6 s, the deformation displacement amplitude simulated by the PMC method is slightly lower than the experimental results. Therefore, within the allowable error range, the simulation results are consistent with the actual physical results.
The original experiment is to verify the multi-physical field coupling problem between water flow, air, and structure. This paper only considers the FSI problem between fluid and structure. The scale of simulation is relatively small, and the accuracy of the SPH solver itself is limited. Therefore, comprehensively excluding the influence of the above factors, it can be concluded that the PMC method proposed in this paper is effective in solving the FSI problem.

\subsection{Dam break flow passes through an elastic sluice gate}
\begin{figure*}[!htb]
	\centering
	\includegraphics[scale=.8]{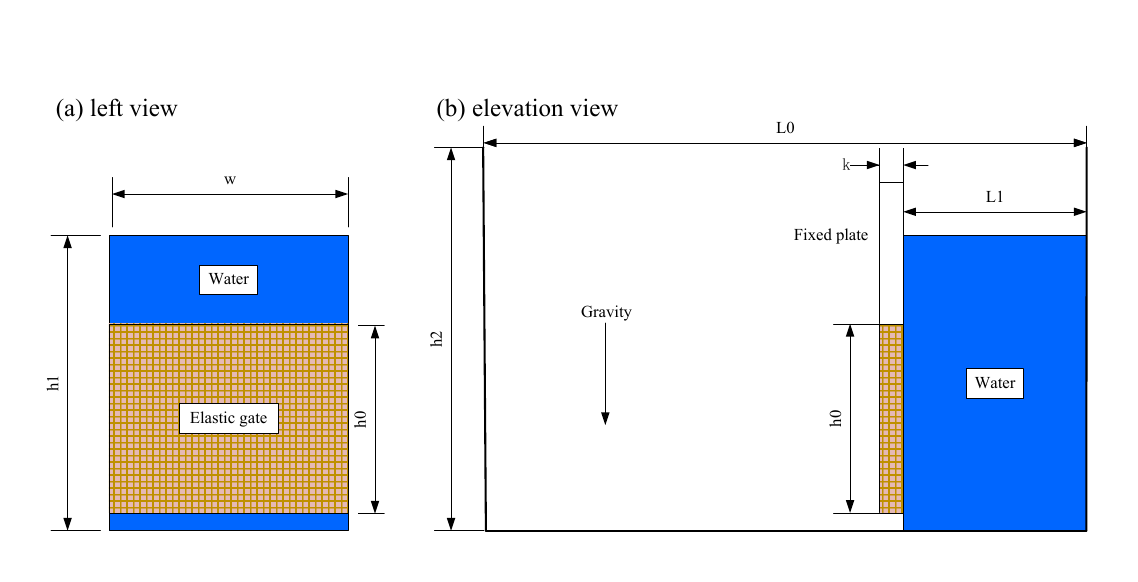}
	\caption{The dimensions of the water tank and the elastic gate, where (a) are the left view and (b) are the elevation view.}
	\label{fig14}
\end{figure*}
\begin{figure*}[!htb]
	\centering
	\includegraphics[scale=1]{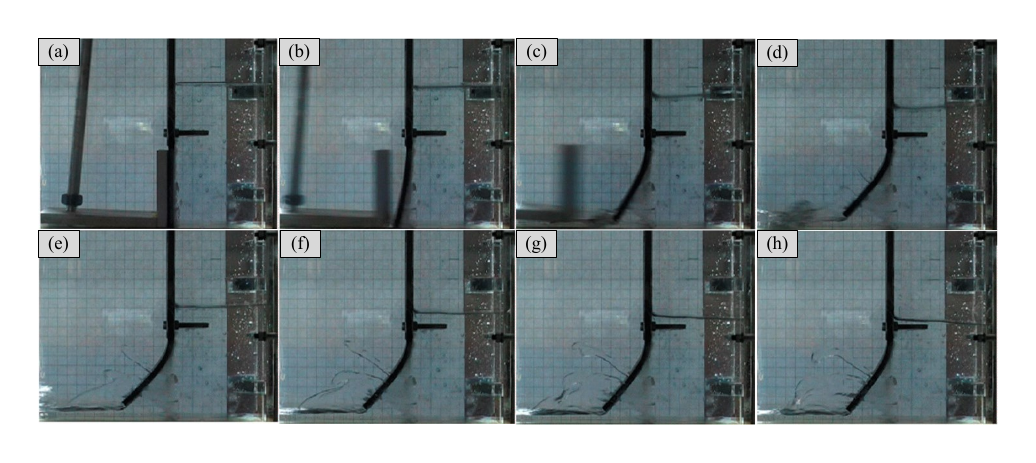}
	\caption{The experimental results of Antoci et al.\cite{antoci2007numerical} from 0 s to 0.32 s, each time interval is 0.04 s.}
	\label{fig15}
\end{figure*}
\begin{figure*}[!htb]
	\centering
	\includegraphics[scale=.37]{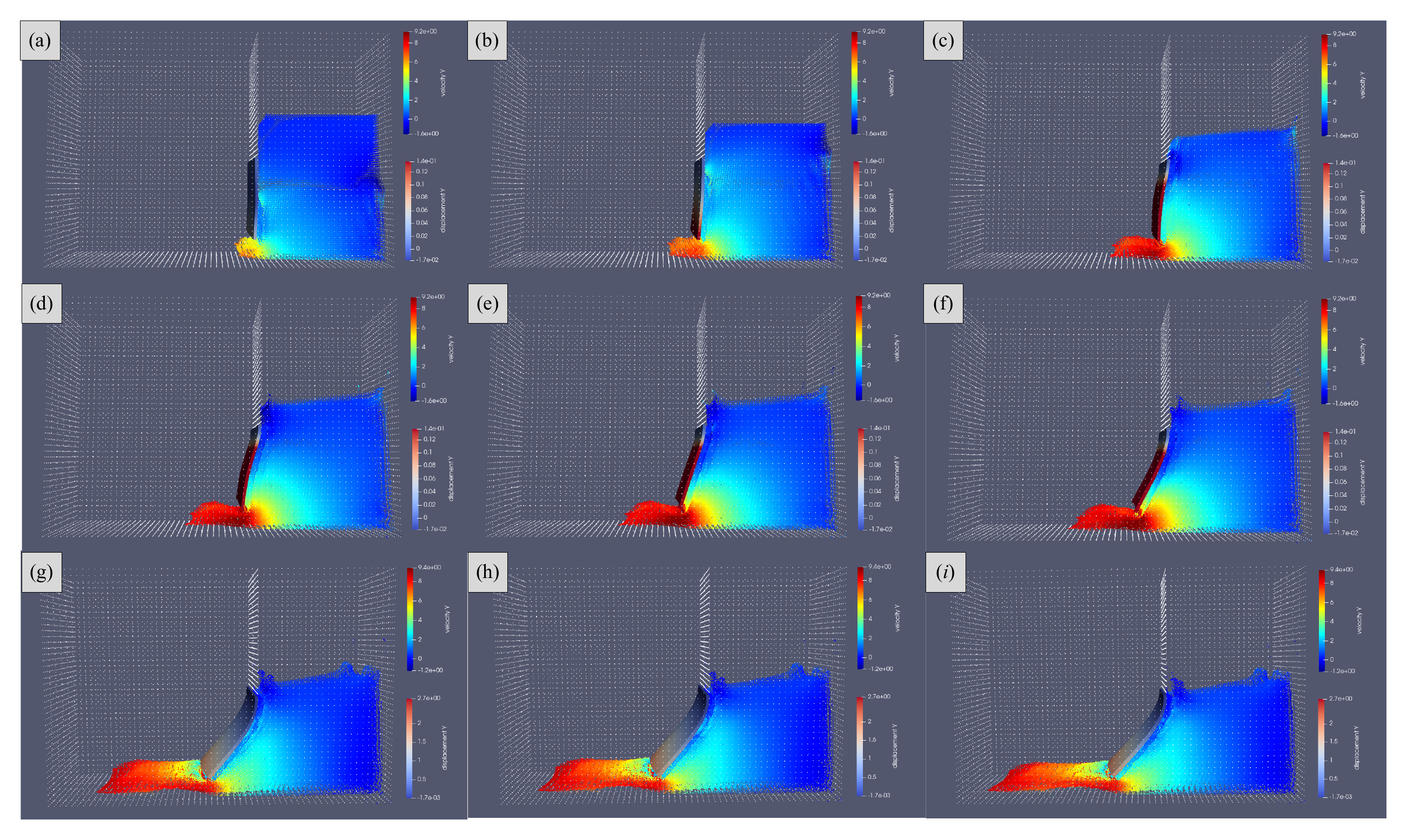}
	\caption{An example using PMC method to simulate water flow through an elastic gate with a time span from 0.4 s to 0.8 s.}
	\label{fig16}
\end{figure*}

In order to further verify the performance of the PMC method, we give the numerical simulation results of water flow through the elastic gate. This example is also a representative benchmark for the FSI problem. Its specific experiment was completed by Antoci \textit{et al.}\cite{antoci2007numerical}. The experimental setup is shown in Fig.\ref{fig14}. The fixed plate and the rubber gate fix the water on the right side of the water tank. The top of the rubber plate is fixed and its bottom is free to move. The whole device is placed in a water tank of size $4\times12\times8$ m. In this experiment, the FEM method is used to disperse the rubber gate, and the SPH method is used to disperse the water flow. The fluid adapter module interpolates data (forces and displacements) through the critical grid and interacts with preCICE.

The experimental dimensions used in this paper are shown in Table \ref{tab03}, where the young's modulus of the elastic plate $E=1.2\times 10^7$ Pa. It determines the deformation of the elastic door under the pressure of water flow. The density of the elastic gate made of rubber is $\rho=1100$ kg/m$^3$, and the Poisson's ratio is $\nu=0.45$. The whole device has a vertical downward gravity field. The value of gravitational acceleration in this paper is $\bm{g}=9.8$ m/s$^2$, the actual results of this experiment are shown in Fig.\ref{fig15}.
\begin{table*}[thp]  
	\centering
	\caption{Geometric parameters of dam break flow passes through an elastic sluice gate (unit: meters).}  
	\label{tab03}
	\begin{tabular}{p{1.8cm}p{1.8cm}p{1.8cm}p{1.8cm}p{1.8cm}p{1.8cm}p{1.8cm}}
		\toprule
		Length of water tank $L_0$ (m) & Height of water tank $h_2$ (m) & Width of water $w$ (m) & Length of water $L_1$ (m) & Height of water $h_1$ (m)& Height of elastic plate $h_0$ (m)& Thickness of elastic plate $k$ (m)\\
		\midrule
		12 & 8 & 4 & 4 & 6 & 3 & 0.1\\
		\bottomrule
	\end{tabular}
\end{table*}

The experiment of  Antoci \textit{et al.}\cite{antoci2007numerical} showed that the whole fluid is in the initial state at t=0 s, the rubber gate is fixed by the baffle, and the water generates pressure on the surrounding wall. With the removal of the baffle, the bottom of the rubber gate began to deform greatly under the pressure of water at t=0.04 s. As the water flows out of the bottom of the rubber gate, the rubber gate is further deformed by the pressure of the water flow, and the state is maintained until the water flows out from the right side. In this procedure, the elastic gate has been under the pressure of water from the right side, so the bottom of the elastic gate deforms to the left.

Fig.\ref{fig16} is the simulation result using the PMC method with a time span from 0.3 s to 0.8 s. Obviously, rubber gate begin to be subjected to the pressure of water flow and has a slight distortion. The bottom of rubber gate is movable, so it will be subjected to the pressure of water flow and show greater distortion at the bottom. As the pressure on the bottom of the right flow increases, the water flows out from its bottom and reaches a relatively stable state. As the fluid particles on the right decrease, the pressure on the bottom of the rubber gate decreases until a new force balance is reached and the elastic plate finally returns to its original position. The simulation results are in good agreement with the experiments of Antoci et al. Because the precision of the SPH solver is limited, the size used in this paper is slightly larger than that of the actual experiments. Therefore, the results of each phase are lagging behind those of the actual experiment. However, the simulation results show that the flow of water and the deformation procedure of the elastic plate is consistent with the experimental results.
\begin{figure*}[!htb]
	\centering
	\includegraphics[scale=.2]{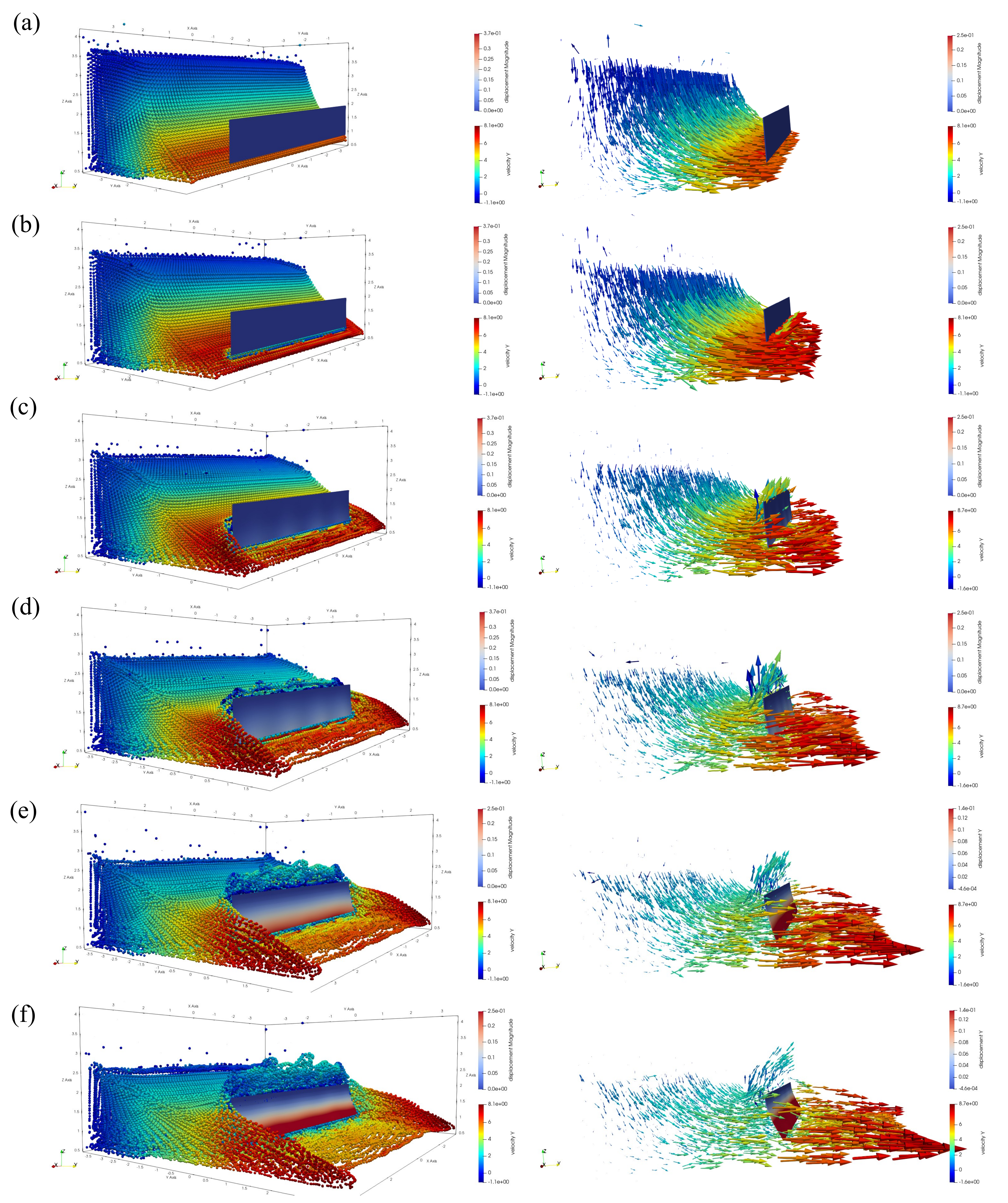}
	\caption{The actual physical time of the simulation process is from 0.5 to 1 second, where each time interval is equal and 0.1 second. The left side is the position and velocity information of the fluid particles (color indicates the velocity), and the right side is the flow field corresponding to the particles.}
	\label{fig17}
\end{figure*}

Fluid-structure interaction (FSI) problem is a classical problem in multi-physical field coupling, which has been successfully applied in many fields  \cite{liu2019smoothed,zhang2021multi,gong2016two}. The fluid-structure interaction phenomenon is also widespread in practice. The study of fluid-structure interaction has very important applications in aerospace \cite{takizawa2020computational,vijayanandh2019design}, water conservancy engineering \cite{wang2017numerical,guo2020fluid,guo2020fluidpp}, safety protection \cite{paik2018one}, petrochemical \cite{lee2018fluid}, oceanographic engineering \cite{wilkes2016numerical}, and biological engineering \cite{samaee2017coupling}. For example, the earthquake disaster area is prone to cause debris flow disasters, which can destroy houses, bridges, and other public facilities and pose a major threat to human life and property. In order to study the impact force of debris flow on the target buildings and establish stable safety protection measures, Dai \textit{et al.} \cite{dai2017sph} established the FSI model of simulated debris flow, in which the debris flow material is simulated as a viscous fluid and solved by SPH, while the check dam is simulated as elastic solid. They solved the governing equations of the two phases and calculated the interaction between them. Finally, they verified the model by simulating the sand flow model test and successfully applied its propagation evolution process to practical problem research. In addition, Zhang \textit{et al.}\cite{zhang2017smoothed} summarized the application of the FSI method in oceanographic engineering in detail, including the capture of phenomena such as water waves, impact, and splash jets. The formation of these phenomena is mainly related to free surface flow.

For most computational fluid dynamics (CFD) solvers based on the Euler framework principle, such as the finite element volume method and the finite difference method, it is a difficult and challenging task to simulate these complex free surface flows. As a Lagrangian and meshless method, smoothed particle hydrodynamics (SPH) has the ability to track different complex boundaries and provides direct support to meet different boundary conditions. Due to the robustness and high accuracy of the SPH method in simulating complex fluid dynamics problems such as free surface boundary, multiphase interface, or material discontinuity, it is very suitable to use the SPH method in fluid simulation of FSI problems. Therefore, another example of dam break impact is given to verify the free surface flow in FSI.

As shown in Fig.\ref{fig17}, it describes the process of free surface flow and solid deformation in the FSI problem. The fluid module is simulated by SPH method and the solid structure is simulated by FEM method. The dam break water flows out from the left side, the top of the square baffle is fixed, and the bottom is free to move. The physical time corresponding to the evolution process from Fig.\ref{fig17} (a) to Fig.\ref{fig17} (f) is from 0.5s to 1s, with a time interval of 0.1s. First, from this evolution procedure, it can be seen that the fluid contacts the solid structure and compresses it, resulting in the deformation of the solid. The water flow at the bottom flows out from the bottom of the baffle. The water flow at the top is blocked by the baffle and generates broken waves. The flow field at the top splashes from above the baffle and flows forward across the baffle. This procedure can capture the actual physical phenomena such as broken waves, the impact of water flow on the baffle and the deformation of the baffle. This example shows that the PMC method proposed in this paper can also effectively simulate a series of physical procedures in FSI problems involving free surface flow.

\section{Conclusion and outlook}\label{section7}

Multi-physics simulation plays an important role in engineering computing. PreCICE, as a multi-physics coupling library, has attracted more and more attention, but so far it only provides an interface to support mesh-based methods. This paper proposes and implements a particle-mesh coupling (PMC) method with preCICE. A critical mesh is designed as an intermediate medium to interpolate and exchange data. This method supports the coupling of the particle method and mesh method to solve FSI problems. Based on the PMC method, we have implemented the adapter supporting the SPH method for the first time, which enables it to connect with the preCICE coupling framework and successfully apply it to solving FSI problems. At the same time, we have given an interpolation method based on the particle and mesh method, which can map the data between fluid and solid with high accuracy, including force and displacement. Finally, we give a verification test of the FSI problem to test the performance of the PMC method. The test results show that the PMC method can effectively solve the FSI problem. Since the computational efficiency of the current version has not been fully optimized, the future work of this paper is to focus on the improvement of the computational efficiency of the PMC method. We will develop PMC programs that can run efficiently on modern supercomputers, and further expand the application scope based on meshless method coupling.

\begin{acknowledgments}
This work was supported by the National Natural Science Foundation of China (Grant No. 61902413 and No. 62002367).
\end{acknowledgments}

\section*{References}
\bibliographystyle{unsrt}
\bibliography{refs}

\end{document}